\documentclass[authoryear]{elsarticle}
\usepackage{comment}
\usepackage{tcolorbox}
\usepackage{amsmath}
\usepackage{amssymb}
\usepackage{amsfonts}
\usepackage{multirow}
\usepackage{booktabs}
\usepackage{tabularray}
\usepackage{wrapfig}
\usepackage[strings]{underscore}
\usepackage{url} 
\usepackage{xurl}
\usepackage{enumitem}
\usepackage{natbib}
\usepackage{xparse}
\usepackage[colorlinks=true,linkcolor=blue,citecolor=blue,urlcolor=blue]{hyperref}
\usepackage{natbib}
\usepackage{xparse}
\usepackage{needspace}
\RenewDocumentCommand{\cite}{m}{%
  \begingroup
    \def\mysep{}%
    \renewcommand*{\do}[1]{%
      \mysep\citeauthor{##1}, \citeyear{##1}%
      \def\mysep{; }%
    }%
    \docsvlist{#1}%
  \endgroup
}

\usepackage{amssymb}
\journal{Elsevier}

\begin{document}

\begin{frontmatter}



\title{Hedging Hydrogen: Planning and Contracting Under Uncertainty for a Green Hydrogen Producer}


\author[label1,label2]{Owen Palmer}
\author[label2]{Hugo Radet}
\author[label1]{Simon Camal}
\author[label3]{Jalal Kazempour}
\author[label1]{Robin Girard}

\affiliation[label1]{organization={Mines Paris - Centre PERSEE},addressline={1 Rue Claude Daunesse},city={Valbonne},postcode={06560},country={France}}
\affiliation[label2]{organization={Verso Energy},addressline={49b Av. Franklin Delano Roosevelt},city={Paris},postcode={75008},country={France}}
\affiliation[label3]{organization={Technical University of Denmark (DTU)},addressline={Anker Engelunds Vej 101},city={Kongens Lyngby},postcode={2800},country={Denmark}}

\begin{abstract}
Green hydrogen production by water electrolysis using renewable electricity is considered essential for decarbonisation of certain sectors of the global economy, however development of the industry is lagging behind expectations due to the perceived financial risk for individual projects. This risk stems from a number of uncertainties, including future hydrogen demand, variable renewable energy sources, and volatile energy market prices.
The interaction of these uncertainties is complex, yet the analysis of hydrogen projects is often carried out using simplified modelling that often omits uncertainty and/or energy hedging practices which are typical for intensive power consumers. In this study, we define a set of planning methods (planning policies) in order to compare the effectiveness of different modelling approaches. We propose a 2-stage market-focused stochastic program to represent a hydrogen producer supplying an industrial customer through a hydrogen offtake contract (a Hydrogen Purchase Agreement, or HPA). The model can be used to obtain equipment sizing decisions, as well as energy hedging decisions using Power Purchase Agreements (PPA's) and power futures. We find that for some HPA contract types, failure to use stochastic modelling can lead to planning decisions that result in 30\% higher production costs during scenario stress-testing for the same project. This could lead to some projects being discarded by developers, incorrectly deemed to be unviable due to cost projections being too high. The results also show the importance of HPA contract volumetric obligations in limiting demand uncertainty.

\end{abstract}

\begin{graphicalabstract}
\includegraphics[scale=0.6]{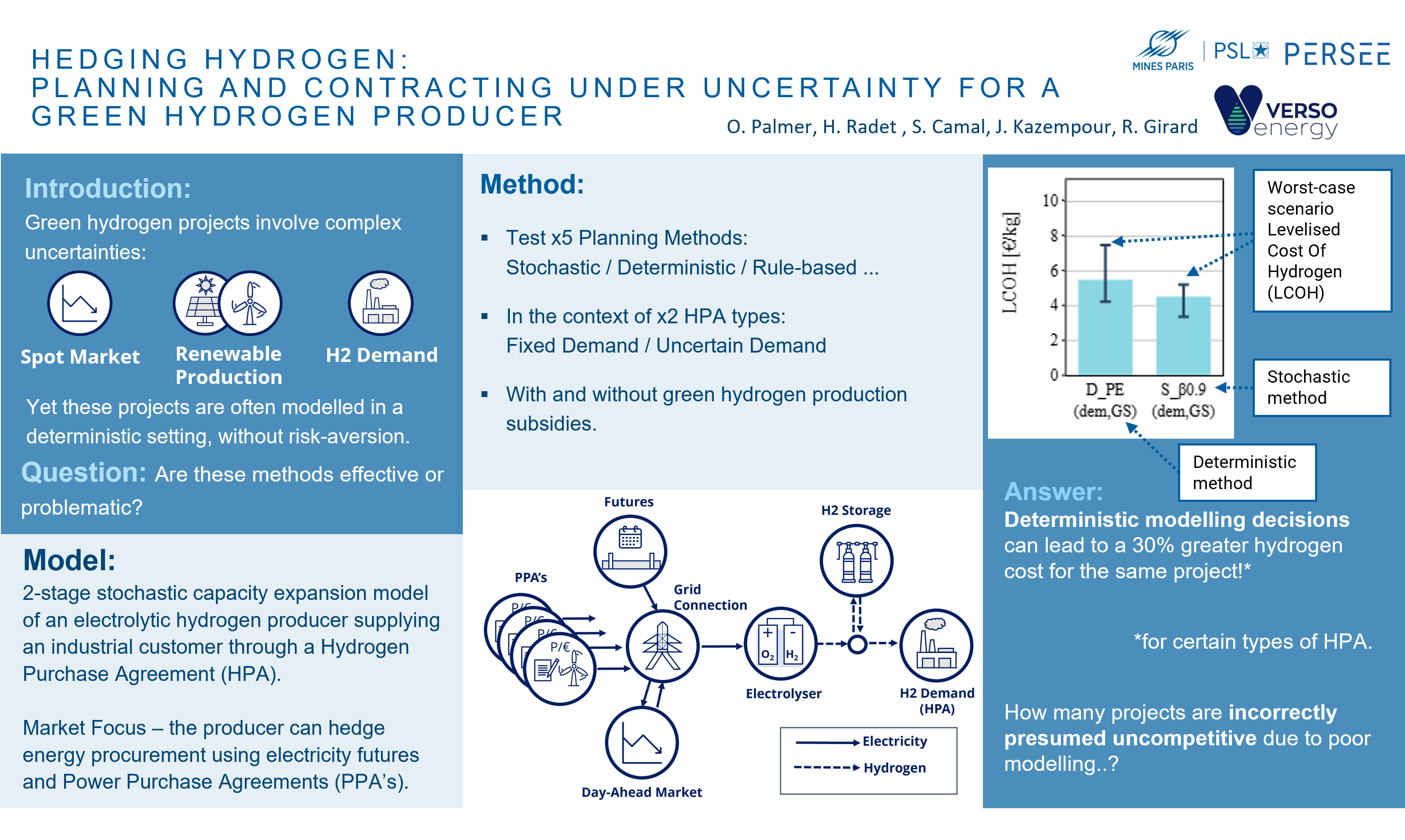}
\end{graphicalabstract}

\begin{highlights}
\item An energy market-focused 2-stage stochastic model of an electrolytic hydrogen plant.

\item Co-optimisation of sizing and energy hedging for Hydrogen Purchase Agreements (HPA).

\item Deterministic planning methods may cause viable projects to be overlooked.

\item Fixed-demand HPAs lower CAPEX, while uncertain demand raises hydrogen sale prices.

\end{highlights}


\end{frontmatter}



\section{Introduction} \label{sec:intro}
\subsection{Context} \label{subsec:intro:context}
As part of the energy transition and the fight against climate change, expanding the use of green hydrogen produced from water and renewable energy through electrolysis is considered to be a promising method of reducing emissions in many hard-to-decarbonise sectors of the economy. Aside from decarbonising current use cases of hydrogen such as fertiliser production and oil refining, green hydrogen may be used in the future to replace fossil fuels in logistics, in the creation of sustainable shipping and aviation fuels, as well as a chemical agent or feedstock in certain industrial processes (\cite{InternationalEnergyAgency2019TheOpportunities}). 

Development of the \textit{low-carbon hydrogen economy} (\cite{InternationalEnergyAgency2019TheOpportunities}) requires investment in both supply and demand. Production necessitates investment in expensive electrolysers, and use of hydrogen on the demand side often requires re-tooling or new equipment. There is also considerable uncertainty in production costs for green hydrogen, as they are closely related to electricity prices which are highly volatile and sensitive to geopolitical shocks.

Attracting the necessary financing to cover these high capital costs can be difficult given the risk-aversion of lenders with respect to the aforementioned risks (\cite{GreenHydrogenOrganisation2022GreenProjects}), as well as the high cost of capital and slow implementation of regulation (\cite{InternationalEnergyAgency2023Global2023}). As a result, green hydrogen project investment is lagging behind the pace and volume required to meet global climate ambitions (\cite{InternationalEnergyAgency2023Global2023,Odenweller2024TheGap}).

Green hydrogen producers require long-term confidence in demand for their product to make their projects \textit{bankable} and attract the necessary financing to cover significant capital costs (\cite{Craen2023FinancingProject,Scholvin2025De-riskingAfrica}). Especially given that markets for green hydrogen do not currently exist, the natural solution for this uncertainty is to find industrial customers willing to engage in bilateral \textit{offtake} contracts, or \textit{Hydrogen Purchase Agreements} (HPA's), wherein a future producer and a future consumer of hydrogen commit to exchanges of hydrogen over an extended period (\cite{Craen2023FinancingProject,Scholvin2025De-riskingAfrica}). 

In order to sign up to an HPA, green hydrogen consumers require long-term confidence in competitive prices for green hydrogen to avoid \textit{“the risk of becoming locked into an expensive and potentially scarce energy carrier”} (\cite{Odenweller2024TheGap}). Potential consumers must be reassured that their hydrogen purchase price will stay competitive with alternatives over the life of the contract even in adverse energy market conditions. 

One missing piece of the puzzle is therefore for potential green hydrogen producers to demonstrate their ability to maintain Levelised Cost of Hydrogen (LCOH) prices within a competitive range for the long-term, despite a strong presence of uncertainty in production costs. Unfortunately, many commonly used quantitative methods for hydrogen project techno-economic analysis in both scientific literature and industry do not fully take into account the impact of this uncertainty. This leaves project developers vulnerable to the following errors: (1) not developing an otherwise viable project because an ineffective planning method deemed it to be unviable, or (2) developing a project with sub-optimal planning decisions that do not provide robustness in the face of future events.

The aim of this paper is therefore to investigate the relative effectiveness of different types of quantitative planning methods (or Planning \textit{Policies} using the terminology of \cite{Powell2022SequentialII}) for decision making in green hydrogen production in the context of an HPA contract.

\begin{figure}
  \centering
    \includegraphics[scale=0.6]{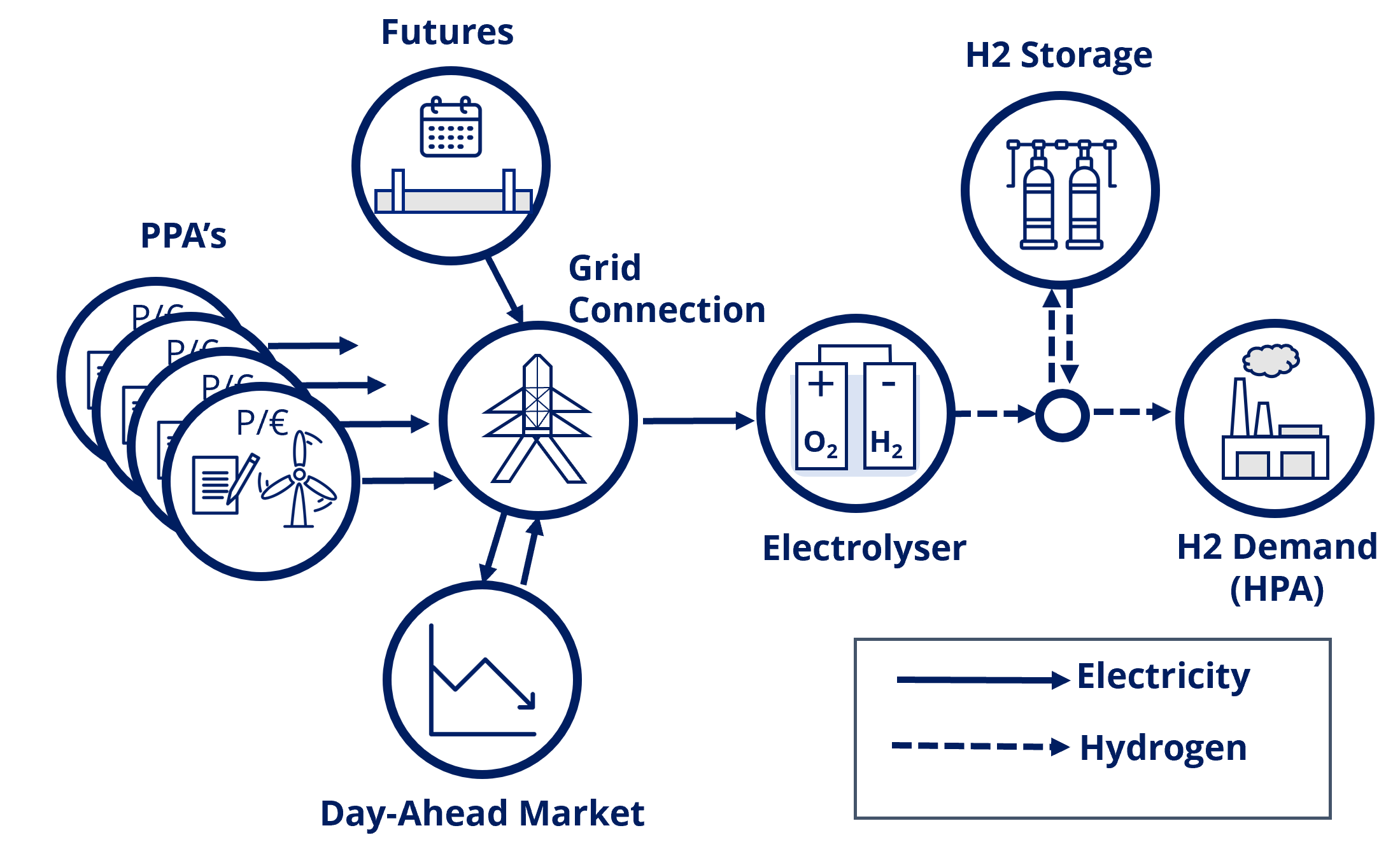}
    \caption{Case study - H2 production model, energy flows diagram.}
    \label{fig:schema}
\end{figure}

With the aim of reducing financial risk, green hydrogen producers have several risk mitigation mechanisms at their disposal which can be separated into the two categories of asset sizing, or energy hedging. Asset sizing decisions may include oversizing of certain flexible assets (such as hydrogen storage) to allow increased utilisation of renewable sources or to ride through market price peaks. Energy hedging decisions relate to the procurement of electricity in advance outside of the spot markets. The most commonly used hedging instruments to date being power futures (where electricity prices are locked in for block volumes with delivery at a future date) (\cite{Sanda2013SelectiveCompanies}), although renewable \textit{Power Purchase Agreements} (PPA's) are becoming more and more common (\cite{Flottmann2024DerivativesMarket}). PPA's are non-standard contracts for which the most typical form, \textit{pay-as-produced}, locks in a set price at which a set proportion of production of a given renewable energy park is bought by the consumer as it is generated (\cite{Mittler2023ReviewAnalysis}).

Additional help may be available in the form of incentives for green hydrogen producers (\cite{Hoogsteyn2025InteractionsHydrogen}), such as subsidy auctions in Europe (\cite{EuropeanCommission2023QuestionsHydrogen}), and a production tax credit in the US (\cite{USDepartmentofEnergy2023AssessingCredit}). Eligibility for these subsidy schemes is to require hourly time-correlation of hydrogen production with renewable sources (\cite{Guillotin2025HydrogenComparison}). As such, the benefit of achieving eligibility to receive these subsidies will influence electricity procurement decisions for individual projects.

The interplay between these sources of uncertainty and potential instruments for financial risk mitigation is far from simple, and yet in both industry and academia approximations are commonly used in order to simplify modelling and decision-making (which are discussed in the literature review to follow). This study builds on literature in three thematic areas: (1) techno-economic analysis of hydrogen electrolysis projects, (2) capacity expansion planning under uncertainty, and (3) electricity hedging under uncertainty.

\subsection{Literature Gaps and Paper Objectives}
\subsubsection{Techno-Economic Analysis of Hydrogen Projects}\label{subsec:intro:lit:analysis_hydrogen_projects}
Many techno-economic feasibility and capacity expansion studies in the hydrogen domain can be found, however the majority are performed using deterministic modelling which inherently ignores uncertainty, and overlooks the risk-averse preferences of investors and project developers. The impact of these simplifications in the context of hydrogen modelling is yet to be quantified in a scientific study. In reviews of modelling tools available for energy system modelling  (\citet{Haugen2023RepresentationReview}), and specifically in the hydrogen domain (\cite{Reulein2025TheModels}), it is noted that “most applied power system tools are deterministic” (\citet{Haugen2023RepresentationReview}). Anecdotal evidence suggests that deterministic modelling is indeed common in industry for project planning and valuation.

Some examples of this in scientific literature include (\cite{AliKhan2021DesigningAustralia,Matute2023Techno-economicPPAs}), in which simulation-based feasibility studies for green hydrogen production via pay-as-produced PPA's is presented, although in a deterministic setting and without optimising capacities. A similarly deterministic model without capacity expansion is proposed in \citet{Koleva2021OptimalMarkets}, in which the authors use iterative solutions of the model to improve electrolyser capacity solutions for different energy procurement configurations. \citet{Brandt2024CostFramework} use a capacity expansion model to study hydrogen cost and emissions sensitivity to electricity prices, however the model used is again deterministic. 

Several studies examine system-level impacts of hydrogen policies and incentives, although also in a deterministic setting (\cite{Ricks2023MinimizingStates,Shirizadeh2023Long-termHydrogen,Zeyen2024TemporalHydrogen,Vargas-Ferrer2025ComplyingChains,Hoogsteyn2025InteractionsHydrogen}). As noted by \citet{Hoogsteyn2025InteractionsHydrogen} in the limitations of their study, the assumption of perfect foresight (which is implicit to deterministic modelling) prevents deterministic models from taking into account risk-mitigating behaviour of investors and project developers.

\subsubsection{Hydrogen Capacity Expansion Planning Under Uncertainty}\label{subsec:intro:lit:capacity_expansion}
There are scientific studies of capacity expansion under uncertainty in the hydrogen domain and the energy domain more broadly (\cite{Roald2022PowerApplications}), however no single study takes into account all the necessary sources of uncertainty, flexible assets, and especially hedging instruments, all of which are required for a practical solution of the problem addressed in this paper. 

For example, the studies by \citet{DeWeerdt2023TheDecisions} and \citet{Fabianek2024ACalifornia} apply real options approaches to hydrogen generator retrofitting and to electrolyser installation respectively, however these studies examine the tipping points between different equipment investment decisions, and how these tipping points change over time --- they do not consider hedging of a contracted demand. 

Similar limitations can be seen in studies using stochastic optimisation approaches for hydrogen system planning. \citet{Coppitters2019Surrogate-assistedUncertainty} propose a stochastic capacity expansion model for sizing of an electrolysis plant with uncertainty in both solar production and equipment CAPEX, although neither demand uncertainty nor grid connection is included. \citet{Li2022CoordinatedUtilization} perform a stochastic optimisation of electrolyser and HVDC capacities under spot price and renewable production uncertainty, though demand uncertainty and hedging are not considered.

Stochastic capacity expansion studies outside the hydrogen domain also provide some but not all necessary components. \citet{Mitra2014OptimalModeling} propose a capacity expansion model for an industrial consumer with uncertain demand and energy prices, however variable renewable production and energy hedging are not included. A price-taker stochastic model for risk-averse generation investment is given in (\cite{Pisciella2016AExpansion}), although it does not include demand and renewable production uncertainty, storage, or hedging instruments. 
\citet{Gal2017FuelMarket} study gas generator capacity decisions, comparing with and without rule-based options hedging. The authors note that increasing price volatility leads to increased equipment capacity investment, however their study does not include risk-aversion, hydrogen production, or storage, and it uses a price-maker equilibrium model that is not appropriate for a single market participant interested in hedging their risk.

\subsubsection{Electricity Hedging Under Uncertainty} \label{subsec:intro:lit:hedging}
There is a rich variety of hedging studies in the energy domain, however examples within the sub-domain of hydrogen production are rare. In addition, energy hedging is typically studied in an operational manner (i.e., with fixed physical asset sizes), and co-optimisation of equipment capacity is usually not included. This ignores the impact of equipment sizing on energy procurement decisions for flexible power-intensive consumers like hydrogen producers.

Energy hedging strategies can be defined as \textit{static} (unchanging as new information is revealed), and \textit{dynamic} (changing over time to adapt to new information, \cite{Dimoski2023DynamicRisks}). Hedging volumes are often expressed in terms of the hedging ratio --- the proportion of the expected future volume which is hedged. In the energy industry, hedging decisions using heuristics such as predefined static hedging ratios are common (\cite{Sanda2013SelectiveCompanies}), which typically prescribe the purchase of hedging instruments incrementally over time to achieve a planned hedging ratio for an expected future volume. 

There are few studies of energy hedging in the context of hydrogen production. A rare example is \citet{Wu2022ElectricityMarkets}, who propose an operational stochastic program including futures hedging for a hydrogen recharging station, although capacity expansion is not included.

The hydrogen production hedging problem is however analogous to other studies relating to hedging for electricity-intensive consumers, although once again capacity expansion is usually not included. This is true for several studies which apply stochastic optimisation in an operational manner in this context (\cite{Carrion2007AConsumers,Conejo2010DecisionMarkets,Zhang2018Long-TermUncertainty}). In an approach similar to that taken in this article, \citet{Gabrielli2022MitigatingOptimization} present a stochastic hedging model with risk aversion to optimise a portfolio of PPA's for demand fulfilment, however this study does not include hydrogen production, demand uncertainty, futures hedging, or storage.

Studies considering the problem of hedging for an electricity retailer also have useful parallels for the hydrogen planning problem studied in this paper, though these studies usually do not incorporate all the necessary elements. For example, (\cite{Boroumand2015HedgingRetailers,Ernstsen2017HedgingCase,Tegner2017Risk-minimisationConsumption}) use similar methodologies to compare static hedging policies using risk metrics like Value at Risk (VAR) and Conditional Value at Risk (CVAR) for electricity retailers with and without means of generation, however these studies do not include capacity expansion, and variable renewable energy is not included.

Other works analyse hedging from a power producer’s perspective (\cite{Pineda2013UsingProducers,Pircalabu2017JointApproach,Tranberg2020ManagingAgreements,Thakur2023PricingContracts}). While these studies provide relevant techniques with respect to pricing and scenario construction, they are essentially the inverse problem of the consumer-side problem addressed in this paper.

Some studies examine only the impact of hedging on profit and loss (P\&L) with respect to an open position on the spot market, however this does not take into account the impact of demand profiles on optimal hedging strategy. One example of this is \citet{Alexander2013TheSpread} which compares the effectiveness of several hedging policies for crude oil products. The authors note that simple hedging strategies can be effective in certain circumstances. In the electricity domain, a similar study is carried out by \citet{Hanly2018TheMarkets}. In addition to the omission of demand profile effects, these P\&L studies also omit renewables, capacity expansion, and storage, and are therefore of limited use for the hydrogen planning problem.

Co-optimisation of capacity expansion planning and procurement strategy is scarce in scientific literature, and the rare examples performing this use system-level equilibrium models which are not suitable for a single actor interested in hedging their risk. One example is the study by \citet{Kazempour2012StrategicMarkets}, in which generation capacity investment and futures trading are optimised. However its agent-based modelling requires assumptions about maximum capacities available to different market players, which is not realistic for a single price-taker with limited information. A different approach is taken by \citet{Shu2023BeyondObligations}, who use iterative solutions of a stochastic equilibrium model to optimise generator investment and different hedging behaviour under different risk-aversion levels and policy obligations. Futures, options, and PPA's are included, however hydrogen production and storage are not. The same model is adjusted by \citet{FarhadBillimoria2024HedgingMarkets} to include the possibility of battery storage, and it is then used to study changes in hedging costs in a progressively renewable grid. However, this iterative method lacks guaranteed optimality and the equilibrium model remains unsuitable for single-actor hedging. 

\subsection{Speculation and No-Resale Formulations in Hedging Problems}\label{sec:intro:speculation}
Some studies in the electricity hedging domain can also be further characterised as using a \textit{no-resale formulation}, the effects of which have not yet been studied in literature (to the authors' knowledge). This modelling concept shall be introduced in the following sections. 

\subsubsection{Risk-Neutral Pricing}
An important modelling decision is how electricity futures should be priced with respect to the spot market. The use of arbitrary or unrelated prices (such as in  \cite{Wu2022ElectricityMarkets} and some studies of \cite{Conejo2010DecisionMarkets}) may inadvertently make the purchase of futures attractive for their resale value alone --- not for alleviating risk. This can lead to excessive \textit{speculation} within the model. The term \textit{speculation} is used in this article to refer to the procurement of energy with the express purpose of reselling that energy on the spot market for profit, rather than using it to fulfil the underlying demand process\footnote{This action is often termed \textit{arbitrage} in energy modelling literature (\cite{Conejo2010DecisionMarkets,Kazempour2012StrategicMarkets,Pineda2013UsingProducers}), however this differs from strict definitions of arbitrage in finance theory (\cite{Delbaen2008TheArbitrage}).}.
A commonly used approach is to adjust either the scenarios, or the hedging instrument prices, so that speculation has zero profitability on average, a state that is often referred to as \textit{risk-neutral} pricing (\cite{Deschatre2021AApplications}). \citet{Fleten2003ConstructingMarkets} and \citet{Dimoski2023DynamicRisks} use risk-neutral day-ahead scenarios produced from known futures prices. In (\cite{Carrion2007AConsumers,Boroumand2015HedgingRetailers,Tranberg2020ManagingAgreements}), the reverse is done: risk-neutral futures prices are calculated based on the expectation of day-ahead market simulations from spot price models. 

\subsubsection{No-Resale Formulations}\label{subsubsec:intro:lit:noresale}
If risk-neutral pricing is not used for a given hedging option, the decisions of a risk-neutral optimisation will almost certainly include an element of speculation. As noted by \citet{Boroumand2015HedgingRetailers}, not using risk-neutral futures prices would mean that the futures contract valuation method will become the main driver for the results. 

In order to avoid this problem without using risk-neutral pricing, a common feature of these aforementioned studies (\cite{Carrion2007AConsumers,Zhang2018Long-TermUncertainty,Kazempour2012StrategicMarkets,Wu2022ElectricityMarkets,Matute2023Techno-economicPPAs,Gabrielli2022MitigatingOptimization}) is the limitation (or more often, complete suppression) of hedging energy resale on the spot market. Modelling formulations that limit hedging energy resale shall be referred to as \textit{no-resale formulations}). These formulations represent a departure from reality, as excess energy can usually be sold on the spot market. Indeed in practice it must, as hedging instrument delivery profiles will not cover exactly the profile of the demand to be hedged, and so systematic buying the shortfall and selling the surplus on the spot market will occur even if 100\% hedged on yearly volume. A thesis proposed in this paper is that no-resale formulations should not be used when deciding between different hedging options.

\subsubsection{Summary and Literature Gaps}
Despite the fact that hydrogen projects involve a complex mix of flexible assets, subsidy structures, numerous sources of uncertainty, and are typically developed by risk-averse agents, to the authors’ knowledge no existing study presents a complete method for jointly optimising equipment investment and electricity hedging for these projects. Furthermore, the impact of common modelling practices such as deterministic modelling, risk neutrality, and no-resale formulations is not known, despite the prevalence of these methods in academia and in industry. 

\subsection{Contributions} \label{subsec:intro:contributions}
A summary of the key contributions of this paper are as follows:
\begin{itemize}
    \item An energy market focused 2-stage stochastic model of an electrolytic hydrogen production plant in the context of an HPA, with co-optimisation of plant capacities, multiple electricity sourcing options, and the possibility of green hydrogen subsidies.
    \item A set of Planning Policies that provide different solutions to the model, incorporating different approaches to:
    \begin{itemize}
        \item Uncertainty modelling (deterministic or stochastic);
        \item Risk aversion (risk-neutral or risk-averse);
        \item Hedging energy resale (resale-enabled or no-resale formulations).
    \end{itemize}
    \item A methodology and a set of metrics for comparing the relative effectiveness (benchmarking) of these Planning Policies in the contexts of HPA contracts with and without demand profile uncertainty, and with and without the presence of green hydrogen subsidies.
    \item A characterisation of uncertainty in renewable production, day-ahead market prices, and demand, in the context of electrolytic hydrogen production.
\end{itemize}

The structure of this paper is as follows: Section \ref{sec:methodology} provides an overview of the methodology used and the general model formulation. In Section \ref{sec:policies}, a set of Planning Policies are proposed for solving the model, and metrics for comparing their effectiveness are defined. Section \ref{sec:case_study} describes the components of a case study project, including uncertainty analysis and elaboration of the mathematical model. Section \ref{sec:results_methodo} presents the results of the Planning Policy and HPA comparisons. In Section \ref{sec:conclusion} the results are summarised and several perspectives on the work are given.

\needspace{5\baselineskip}  % 
\section{Methodology} \label{sec:methodology}
\subsection{Overview} \label{subsec:methodology:process}
The objective of the project developer is to provide as low as possible hydrogen prices on average across the life of the project, whilst limiting the severity of hydrogen price increases that may be possible in adverse conditions in the future. In order to achieve this, the project developer seeks to obtain planning decisions (henceforth referred to as \textit{design decisions}) that manage production cost risks, whilst minimising total cost. These design decisions include equipment capacity sizing (for the electrolyser, hydrogen storage, and grid connection) as well as energy hedging decisions (such as committing to PPA's and buying electricity on the futures market).

These design decisions are obtained using a Planning Policy in the context of a specific HPA. 

A Planning Policy defines a specific method of decision-making. They may take any number of forms. They may be rule-based or optimisation-based, and may have different objective functions (e.g., minimise cost, maximise electrolyser load factor, minimise spot market exchanges, etc.). The main focus of this study is to compare, with a consistent objective function, Planning Policies with different approaches to risk (e.g., risk-averse or risk-neutral), and different approaches to uncertainty (deterministic or stochastic). 

HPA's in practice may have different conditions that are designed to reduce demand uncertainty, including but not limited to minimum and maximum yearly, monthly, daily, and hourly volumes. In this study, an HPA with fixed demand corresponds to a contract where the consumer is bound to follow precisely their original profile (i.e., all volumes down to hourly volumes are fixed, and no randomisation is included). Uncertainty can then be added by relaxing some of these rules, for example, by allowing monthly volumes to vary. In this case, monthly volumes are then generated with randomisation across different scenarios, all the while ensuring that the conditions for a fixed annual volume and maximum hourly volumes are respected. The conditions of the HPA thus define an envelope in which demand uncertainty is permitted. 

A clear distinction is made between in-sample scenarios (scenarios with which the decisions of a particular policy are \textit{made}), and out-of-sample scenarios (scenarios on which the decisions of a particular policy are \textit{tested}) (\cite{King2012ModelingProgramming}). By only evaluating design decisions on their ability to handle new realisations of uncertainty that they have not seen before, we can have confidence that our solutions are not “fluke” design decisions which have over-fit the in-sample set.

The methodology followed for obtaining and then testing the design decisions of a single policy within the context of a single HPA is shown in Figure \ref{fig:methodo}, and outlined as follows:

\begin{figure}
  \centering
    \includegraphics[scale=0.7]{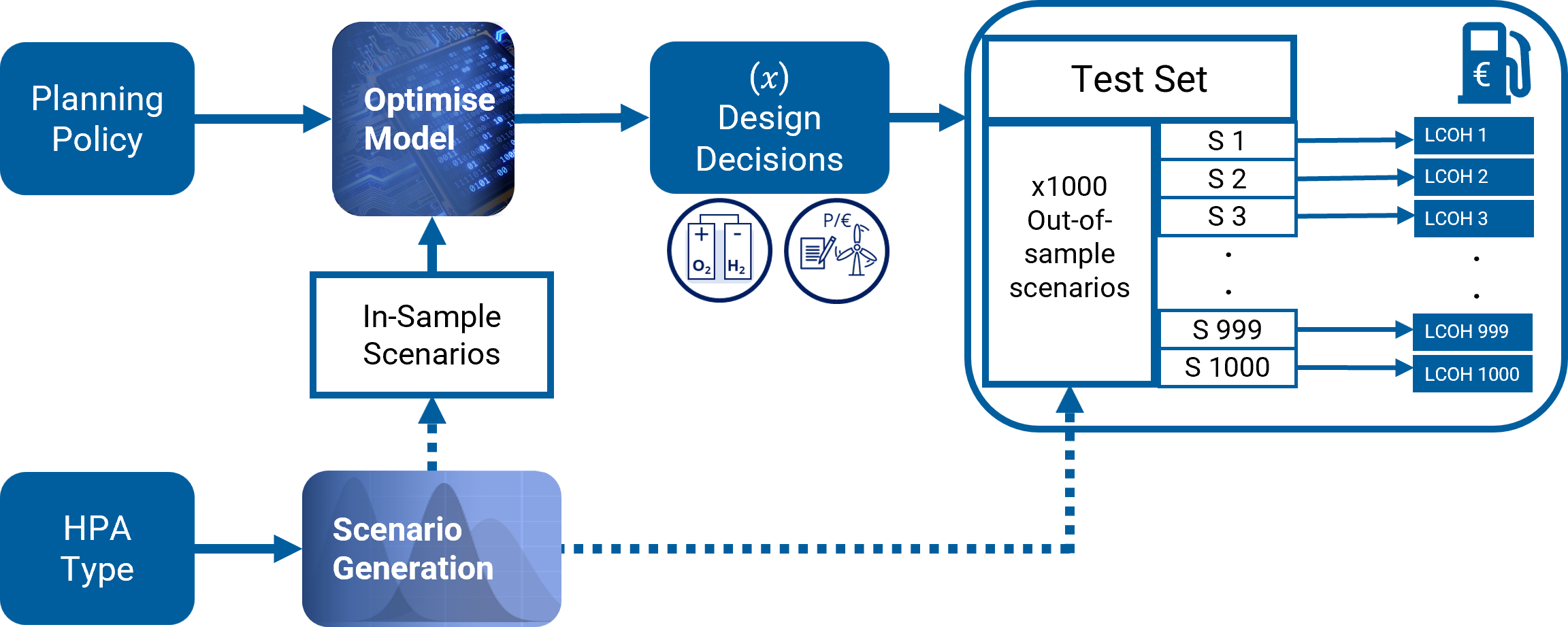}
    \caption{Optimisation and test methodology depicted for a given Planning Policy and a given Hydrogen Purchase Agreement (HPA).}
    \label{fig:methodo}
\end{figure}

\begin{enumerate}
    \item Define the HPA contract type (e.g., fixed demand or uncertain demand, subsidised or non-subsidised) and the uncertainties it allows, to be used to construct the in-sample and out-of-sample scenarios.
    \item Define the Planning Policy to be used to solve the model (see Section \ref{subsec:policies:defintions}).
    \item Use the Planning Policy and the in-sample scenario set to optimise the model and obtain the design decisions of the solution.
    \item Simulate the design decisions on each scenario in the out-of-sample test set and compute its performance (see Section \ref{subsec:methodology:performance_metrics}).
    \item Repeat Steps 2-4 for each Planning Policy to be studied.
    \item Compare the set of performance results for each Policy using the solution comparison metrics (see Section \ref{subsec:methodology:evaluation_metrics}).
\end{enumerate}

In this study, all planning policies were first compared for a Fixed Demand HPA contract. The two best performing Planning Policies were then tested for three more HPA types: a Fixed Demand HPA contract with green subsidies available, and Uncertain Demand HPA contracts with and without green subsidies available.

\subsection{Objective Function} \label{subsec:methodology:obj_func}
The model used is organised as a 2-stage stochastic linear program, with the design decisions $x$ as first-stage decisions. The second-stage ‘operational’ decisions $u_{s,h}$ dictate the dispatching of the assets in one-hour timesteps $h$ for scenario $s$. Operational states such as storage state-of-charge are denoted by $z_{s,h}$. Power values are always given in MW with costs in €/MW, while energy values are given in MWh with costs in €/MWh.

The objective function minimises total annualised cost, in which $J^{d}$ represents the annualised overnight costs of the design decisions, and $J^{o}_s$ represents the net operational costs incurred for scenario $s$ due to spot market exchanges, PPA purchasing, demand curtailment penalties, and network charges, minus the green subsidies received.

The risk aversion factor $\beta$ is used to toggle the weighting between the expected value of the operational cost $J^{o}_{s}$, and the expected value of the worst-case scenario, which is obtained by using the Conditional Value at Risk (CVaR) risk measure (\cite{Rockafellar2000OptimizationValue-at-risk})  with $\alpha=0.99$ (additional CVAR formulation details provided in the Supplementary Material).

\begin{align} \label{eq:obj_fun}
\min_{x, u} \ \underbrace{J^d(x,c,b)}_{\text{Design Cost}}  +   (1-\beta) & E[\underbrace{J^{o}_{s}(x,u_s,w_s,g)}_{ \text{Operational Cost (all)}} ]   \\ \nonumber
+\quad(\beta) & \underbrace{CVaR_{\alpha}[J^{o}_{s}(x,u_s,w_s,g)]}_{ \text{Operational Cost (Worst-Case)}} \\
& \text{s.t.} \nonumber \\
& x \in \mathbb{X}, \quad  u_{s,h} \in \mathbb{U}_{s,h}, \quad z_{s,h+1} \in \mathbb{Z}_{s,h+1},\\ \noalign{\vskip12pt}
& \forall \quad h \in [1..H], \quad s \in [1..S].  \nonumber
\end{align}

The objective function formulation shown in Equation (\ref{eq:obj_fun}) ensures that bias will not occur between operational costs and design costs, and it is easily adapted to provide a number of common planning policies by changing the risk aversion parameter $\beta$. The feasible spaces of the design decisions, operational decisions, and operational states are denoted by $\mathbb{X}$, $\mathbb{U}_{s,h}$, and $\mathbb{Z}_{s,h+1}$, respectively. These spaces are defined in detail for the case study in Section \ref{subsec:case_study:model}.

The same model and objective function are used for optimisation using each planning policy, and for testing of each policy's solution on the out-of-sample set (although for testing the design decisions are fixed).

\section{Planning Policies} \label{sec:policies}
\subsection{Overview} \label{subsec:policies:overview}
The following planning policies emulate common methods used in industry and academic literature and can each be used to provide a solution to the model (i.e., a set of design decisions). They are formally defined in the following section and summarised in Table \ref{tab:policy_details}. The prefixes $D$ and $S$ indicate deterministic and stochastic policies, respectively. The stochastic policies are differentiated by the risk aversion factor $\beta$ used: $\beta0$ for risk-neutral, and $\beta0.9$ for risk-averse. Policies which additionally use a no-resale formulation are identified by the suffix $(NR)$.\\

\begin{table}
\scriptsize
\centering
\caption{Planning policy characteristics.}
\label{tab:policy_details}
\begin{tblr}{
  column{3} = {c},
  hline{1,7} = {-}{0.08em},
  hline{2} = {-}{},
}
\textbf{Policy} & \textbf{Type}            & {\textbf{Risk}\\\textbf{Appetite}} & {($\beta$)} & {\textbf{No-Resale (NR)}\\\textbf{Formulation}} \\
1. $D\_EVP(NR)$            & Deterministic                                        & Neutral                            & N/A                                                                & Yes                                          \\
2. $D\_{PE}$            & Deterministic/Rule-based & Averse                             & N/A                                                                & No                                           \\
3. $S\_{\beta0}(NR)$            & Stochastic                                      & Neutral                            & 0                                                                & Yes                                          \\
4. $S\_\beta0.9$             & Stochastic                               & Averse                             & 0.9                                                              & No                                           \\
5. $S\_\beta0.9(NR)$             & Stochastic                                  & Averse                             & 0.9                                                              & Yes                                          
\end{tblr}
\end{table}

\subsection{Definitions} \label{subsec:policies:defintions}
\textit{1. Expected Value Problem Policy: $D\_EVP(NR)$} \\
The Expected Value Problem policy is the reference case deterministic risk-neutral policy, with co-optimised hedging and equipment capacity decisions. It provides a deterministic solution of the model with one input scenario representing an average case for all uncertainties. The no-resale formulation is applied in order to avoid speculative results (Equation (\ref{eq:op_bounds_last}) is active, as defined in Section \ref{subsec:model:no_arbitrage}). This planning policy is inspired by the prevalence of deterministic modelling in feasibility studies in academia and in industry, as discussed in Section \ref{subsec:intro:lit:analysis_hydrogen_projects}.\\

\textit{2. Pessimistic Expert Policy: $D\_{PE}$} \\
This policy resembles the previous Expected Value Problem policy, although without equipment and hedging co-optimisation. For the hedging decisions, a simple exogenous risk-averse heuristic is used to hedge against high spot market prices: the total yearly hydrogen demand (electrical equivalent), is procured using the least-cost PPA options. It is assumed the expert believes diversity to be important, and so acquires 50\% of the required yearly production from a low-cost solar PPA and the remaining 50\% from a low-cost wind PPA. This calculation is done using the average capacity factor of each park. The remaining design decisions of equipment sizing are then made in a deterministic setting with a single expected value scenario, as for the Expected Value Problem (for uncertain demand this means taking the hour-by-hour average of all in-sample demand scenarios). This policy is inspired by separately optimised hedging decisions and the use of pre-defined hedging ratios, which are in industry as discussed in Section \ref{subsec:intro:lit:hedging}.\\

\textit{3. Stochastic Risk-Neutral Policy: $S\_{\beta0}(NR)$} \\
This policy uses a risk-neutral 2-stage stochastic formulation, with co-optimised hedging and equipment capacity decisions. Operational costs of each in-sample scenario have equal weighting (risk aversion weighting $\beta$ in Equation (\ref{eq:obj_fun}) is set to 0). The no-resale formulation is applied (Equation (\ref{eq:op_bounds_last}) is active). This policy is inspired by risk-neutral stochastic modelling which is typical in capacity expansion models in energy modelling literature as discussed in Section \ref{subsec:intro:lit:capacity_expansion}.\\

\textit{4. Stochastic Risk-Averse Policy: $S\_\beta0.9$} \\
This policy performs a co-optimised risk-averse 2-stage stochastic optimisation, where the expected value of the operational costs and the operational cost of the worst-case scenario are given weightings of 10\% and 90\% respectively (risk aversion weighting $\beta$ in Equation (\ref{eq:obj_fun}) is set to 0.9). This policy is inspired by risk-averse modelling in literature with respect to capacity expansion and energy hedging literature as described in Sections \ref{subsec:intro:lit:capacity_expansion} and \ref{subsec:intro:lit:hedging}. 

While adjustment of either $\beta$ or the CVAR $\alpha$ can achieve similar risk-averse solutions, the main consideration is the use of settings that provide a high level of risk aversion without being fully robust (as this ignores many scenarios in the in-sample set and tends to overfit the worst-case scenario).\\

\textit{5. Stochastic Risk-Averse No-Resale Policy: $S\_\beta0.9(NR)$} \\
This policy performs again a co-optimised stochastic optimisation with the same risk-aversion weighting $\beta$ as for the preceding policy, but with the no-resale formulation applied (Equation (\ref{eq:op_bounds_last}) is active). This policy is included to identify the impact of no-resale formulations as used in energy modelling literature as described in Section \ref{subsubsec:intro:lit:noresale}.\\

\subsection{Solution Notation} \label{subsec:policies:notation}
To distinguish different solutions obtained using the same planning policy, the suffix (GS) is added for Green Subsidies, and (dem) is added for demand uncertainty. 

\subsection{Performance Measures} \label{subsec:methodology:performance_metrics}
The performance of a set of design decisions on a single scenario is measured by the LCOH for that scenario. The LCOH is defined as follows:
\begin{align} \label{eq:lcoh}
LCOH_s = \frac{J^{d} + J^{o}_s}{M^{h2} \cdot \sum^H_{h=1} w^{hd,in}_{s,h}} \quad \text{[€/kg H2],}
\end{align}
where $w^{hd,in}_{s,h}$ represents the hydrogen demand in $MWh$ for scenario $s$ at hour $h$. The mass factor of hydrogen in $kg/MWh$ is denoted by $M^{h2}$, and $H$ is the number of one-hour time steps in the scenario.

The overall performance of a solution obtained by a particular planning policy is then assessed based on its distribution of LCOH results for all scenarios in the out-of-sample test set. The final measures used to assess the performance of a single Policy are: (a) the mean LCOH across all test scenarios, and (b) the worst-case (highest) LCOH for a single scenario in the test set. Worst-case LCOH is an important metric for risk-averse investors as part of scenario stress-testing for project appraisal.

\subsection{Policy Comparison Metrics} \label{subsec:methodology:evaluation_metrics}
The classic metric, the \textit{Value of the Stochastic Solution}, or VSS, (\cite{Birge2010IntroductionProgramming}) is defined as the difference in performance between a deterministic model, and its stochastic equivalent (typically on expected value). In this study, this metric is applied to both expected value and worst-case LCOH performance, and three other similar metrics are proposed: 

The \textit{Value of the Risk-Averse Solution} (VRAS) is defined as the proportional difference in performance of a risk-neutral solution with its risk-averse equivalent.

The \textit{Value of the Resale-Enabled Solution} (VRES) is defined as the proportional difference in performance between solutions optimised using a no-resale formulation, and solutions optimised with unlimited possibility for resale (resale-enabled). 

The \textit{Cost of Demand Uncertainty} (CDU), is defined as the proportional difference in performance between solutions obtained using the same planning policy, but for a Fixed Demand HPA contract, and an Uncertain Demand HPA contract.

\section{Case Study} \label{sec:case_study}
\subsection{Case Study Selection and Applicability}
In this section, a case study of a hydrogen project to supply hydrogen under an HPA in metropolitan France is developed. Specific profiles are used for demand and PPA production, however these may be replaced with any desired profile, making the methodology adjustable to any location with an electricity spot market and hedging instruments comparable to the European market. 

In this case study all design decisions are assumed to be continuous within their allowed range, including both hedging decisions and equipment capacity decisions. Design decisions could be modelled as mixed integer decisions, which is useful for taking into account incremental equipment sizes or fixed PPA contract capacities. However the objective of this study is to identify the inherent dynamics in the system without imposing exogenous limitations. This approach provides useful information at the initial feasibility and conception stage of project development. When progressing to the detailed planning stage, it would then be appropriate to re-optimise the same model proposed in this paper using mixed integer decisions once equipment suppliers and actual available PPA sizes and costs have been identified. 

Previous studies have indicated that the inclusion of lithium-ion batteries is not currently cost-effective for the sole purpose of time-shifting electrical demand of the electrolyser (\cite{Superchi2023Techno-economicSteelmaking,Palmer2023Risk-conscious.}), however the inclusion of value stacking with other services may change this result. Considering additional revenue sources is beyond the scope of this study.

\subsection{Model Components and Uncertainty} \label{subsec:case_study:components_uncertainty}

The model shown in Figure \ref{fig:schema} represents a system configuration for production of green hydrogen to meet demand from a single industrial customer. Electricity used to produce the hydrogen may be procured from the day-ahead market, futures contracts, or from a selection of PPA's.

\subsubsection{Hydrogen Demand} \label{subsubsec:case_study:components_uncertainty:demand}
A hypothetical demand for hydrogen was considered. The HPA contract in all cases specifies an exact annual consumption of 540 t (18 GWh H2), and a maximum hourly demand of 210 kg (6.3 MWh H2). Four different categories of hydrogen projects were considered: 

\begin{itemize}
    \item Fixed Demand HPA;
    \item Fixed Demand HPA -- Subsidised;
    \item Uncertain Demand HPA;
    \item Uncertain Demand HPA -- Subsidised.
\end{itemize}

\begin{wrapfigure}{r}{8cm}
  \centering
    \includegraphics[scale=0.65]{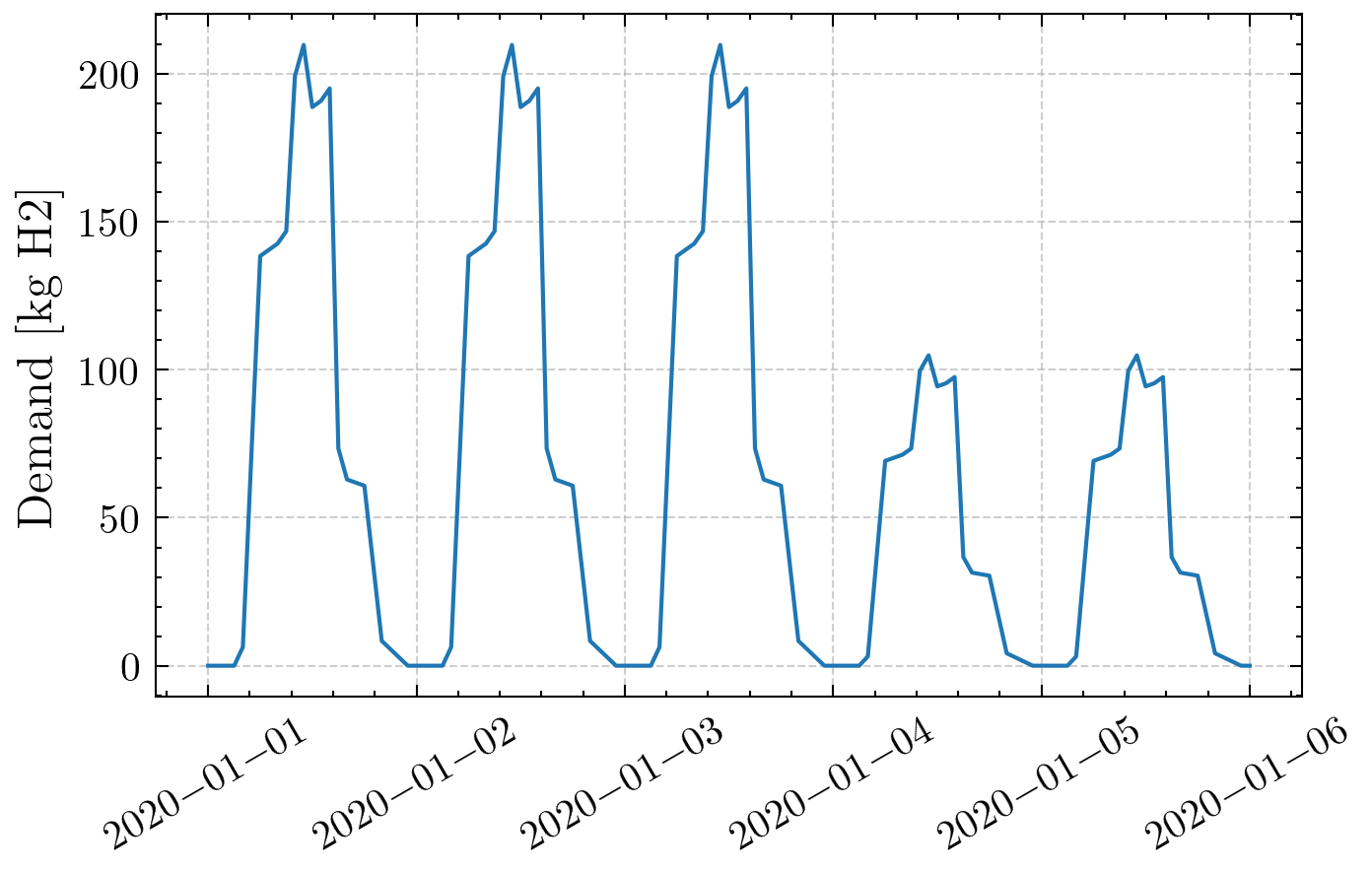}
    \caption{Hydrogen demand Wednesday-Sunday profile for the Fixed Demand HPA.}
    \label{fig:demand}
\end{wrapfigure}

For the case of the Fixed Demand HPA contract, the hydrogen customer does not have the right to vary their future consumption profile. They must obey the original weekday/weekend cyclic profile (see Figure \ref{fig:demand}) with no variations week-to-week or month-to-month. 

\begin{wrapfigure}{l}{8cm}
  \centering
    \includegraphics[scale=0.65]{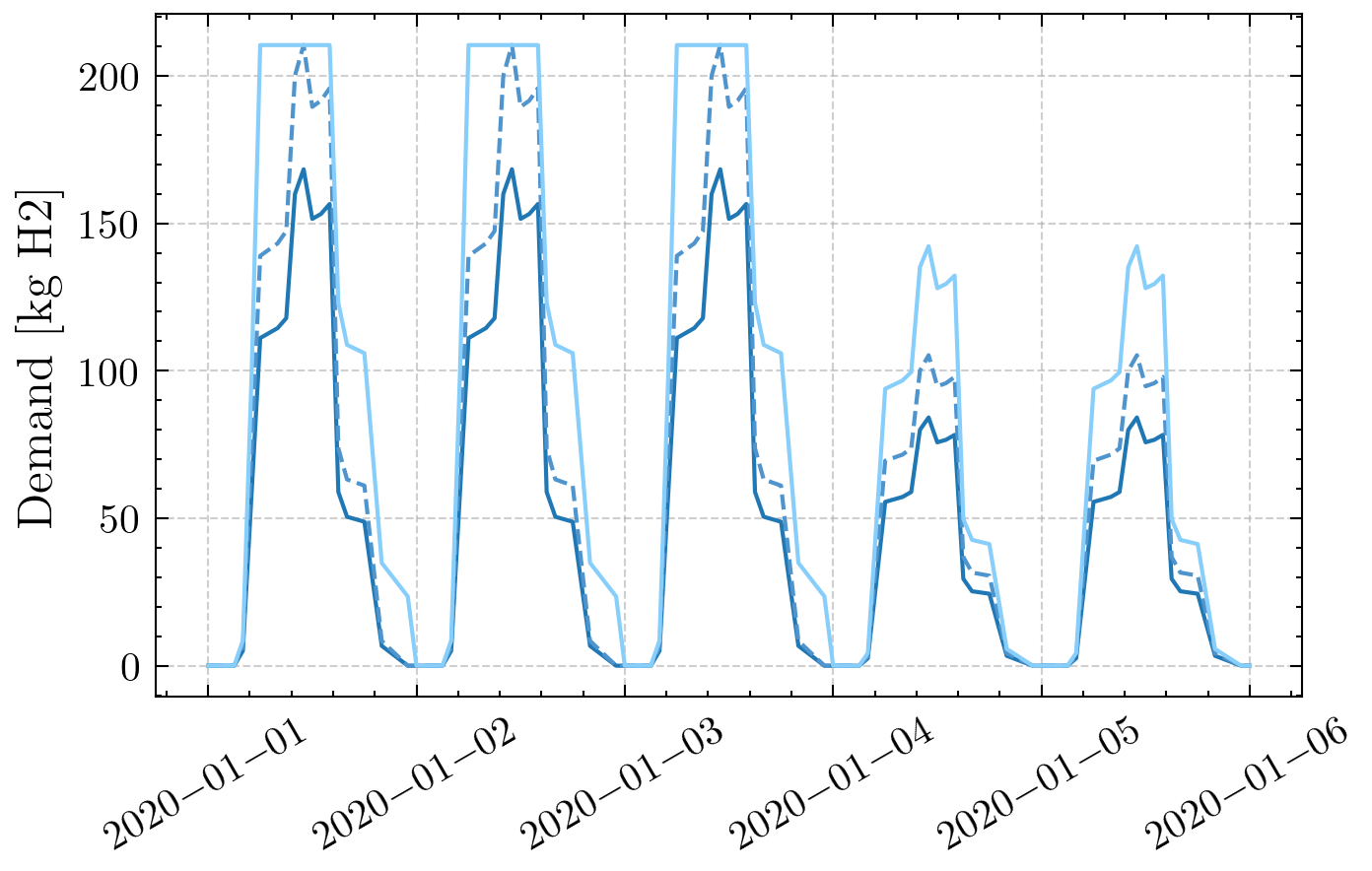}
    \caption{Hydrogen demand Wednesday-Sunday profiles for x3 example scenarios of the Uncertain Demand HPA.}
    \label{fig:demand_scens}
\end{wrapfigure}

For the Uncertain Demand HPA contracts, the consumer will have some flexibility to vary their future consumption profile (i.e., in response to changing demand for their product). It is assumed that the contract dictates the same annual volume and maximum hourly demand as for the Fixed Demand contract, however the customer will have the right to vary their daily, weekly, and monthly volumes. In scenario generation, these volumes are allowed to vary by +/- 50\% with respect to the standard demand scenario used for the strict contract. Three example scenarios are shown in Figure \ref{fig:demand_scens}, and a comparative illustration of their associated month-to-month volume changes is provided in the Supplementary Material.

Another important condition to be defined in the HPA is a monetary penalty for unfulfilled (curtailed) demand. In our case study, for non-delivery of the required demand for a given hour, a penalty of €10,000 /MWh is used in in-sample optimisation to strongly discourage planned use of this option, whilst for testing, curtailment is penalised at €1,000 /MWh, which is commonly used as the electrical curtailment cost in other studies (\cite{Jodry2023IndustrialOptimisation}). 

\subsubsection{Hydrogen Equipment} \label{subsubsec:case_study:components_uncertainty:h2_equipment}
Hydrogen is produced by an electrolyser, followed directly by a compression stage. It may be fed to the customer as it is produced, or be stored and delivered to the customer at a later time. This model uses a simple linear efficiency to model the combined electrolyser + compressor component. This is an approximation of actual electrolysis dynamics which are non-linear, and include discontinuities due to minimum operating levels and distinct operating states (\cite{Zheng2022OptimalDynamics,Baumhof2023OptimizationNecessary}). The impact of including an upper-range minimum operating point of 25\% (\cite{Baumhof2023OptimizationNecessary}) was evaluated using the Expected Value Problem policy and found to result in no significant change in design decisions, and an increase in final cost of less than 0.1\%. It was therefore considered that while modelling these operating states would be preferable, omitting them is not likely to strongly impact the results in this study. 

Neither water flow limits nor water purchasing costs are taken into account. 

\subsubsection{Day-Ahead Market} \label{subsubsec:case_study:components_uncertainty:day_ahead}
Electricity price data for the French day-ahead (spot) market from 2013 to 2022 were obtained (\cite{EEX2022EuropeanEEX,ENTSOe2022ENTSOePlatform}). These scenarios represent a range of average annual prices from 30 €/MWh to 275 €/MWh, with a mean of ~70 €/MWh. The historical distributions of other metrics were also obtained, including average daily spread, and seasonal price drift.

Day-ahead market scenarios were constructed using the historical day-ahead price curves as a base, and then manipulating them using random sampling with respect to historical distributions of these metrics. Additional details of the scenario generation process are provided in Section \ref{subsubsec:case_study:components_uncertainty:scenario_gen} and in the Supplementary Material. 

\subsubsection{Electricity Futures} \label{subsubsec:case_study:components_uncertainty:futures}
The first set of instruments that can be used to hedge against high day-ahead prices is energy futures. The products modelled in this study are those with the longest maturity dates in the European market: calendar year (\textit{CAL}) and quarterly products (\cite{EEX2022EuropeanEEX}). The European Energy Exchange (EEX) also offers two categories of delivery profiles: Baseload Futures, for which the quantity of energy bought is delivered in a constant band for every hour across the delivery period; and Peakload Futures, which provides constant power only during the hours 8am to 8pm. 

In the context of long-term decision making, even these long-maturity futures are only available up to three years in advance, which is clearly not sufficient for informing investment decisions for project lifetimes of up to 25 years. As such, in the context of this study, design decisions that select these products should be interpreted as a strategy to systematically and progressively hedge consumption during the selected period as soon as these products become available for every year of the project life. They are therefore \textit{static global} hedging decisions (\cite{Dimoski2023DynamicRisks}).

Risk-neutral futures prices were calculated as the average of the day-ahead prices for the given delivery period in the relevant scenarios. Peakload futures prices are similarly calculated as the expectation of prices within peak hours within the given period. In practice power futures prices involve an uncertain and time-varying risk premium with respect to day-ahead prices (\cite{Benth2013AnMarkets,Oliveira2021AnalysisAversion,Secomandi2022QuadraticOperations}), however modelling this premium requires making assumptions about market risk factors and their correlation to other variables that are likely to be unreliable and of limited value (\cite{Tranberg2020ManagingAgreements,Dimoski2023DynamicRisks}). Risk-neutral pricing is thus considered to be the most reliable pricing technique for electricity futures for the purposes of this article.

\subsubsection{Power Purchase Agreements (PPA's)}\label{subsubsec:case_study:components_uncertainty:ppas}
The second group of energy hedging products considered is a set of renewable PPA's. A selection of 9 PPA options was included in the model, representing 4 solar photovoltaic (PV) PPA's and 5 wind PPA's in different locations in metropolitan France. Data were obtained for 10 years 2013 to 2022 using  Renewables Ninja (\cite{Pfenninger2016Long-termData,Staffell2016UsingOutput}).

One method for PPA pricing would consider the future value of energy resale on the spot market, bounded below by the Levelised Cost of Energy (LCOE) of the associated generator (\cite{Qorbanian2024ValuationProcurement}). If these two factors were sufficient for predicting PPA prices in the market, PPA prices could be modelled in a risk-neutral manner similarly to electricity futures (possibly including additional risk premiums and option-inspired structures such as in \cite{Tranberg2020ManagingAgreements} and \cite{Pombo-Romero2024AssessingPPAs}). 

In practice, PPA's are non-standard and illiquid, and so PPA prices are highly sensitive to difficult-to-model additional factors including but not limited to counterparty risk, negotiating strength and market power, and green regulations (\cite{Mili2025GreenAgreements}). PPA prices that are above and below their expected fair market value are possible (i.e., they cannot be assumed to be risk-neutral with respect to the spot market). 
\begin{wrapfigure}{l}{8cm}
  \centering
    \includegraphics[scale=0.4]{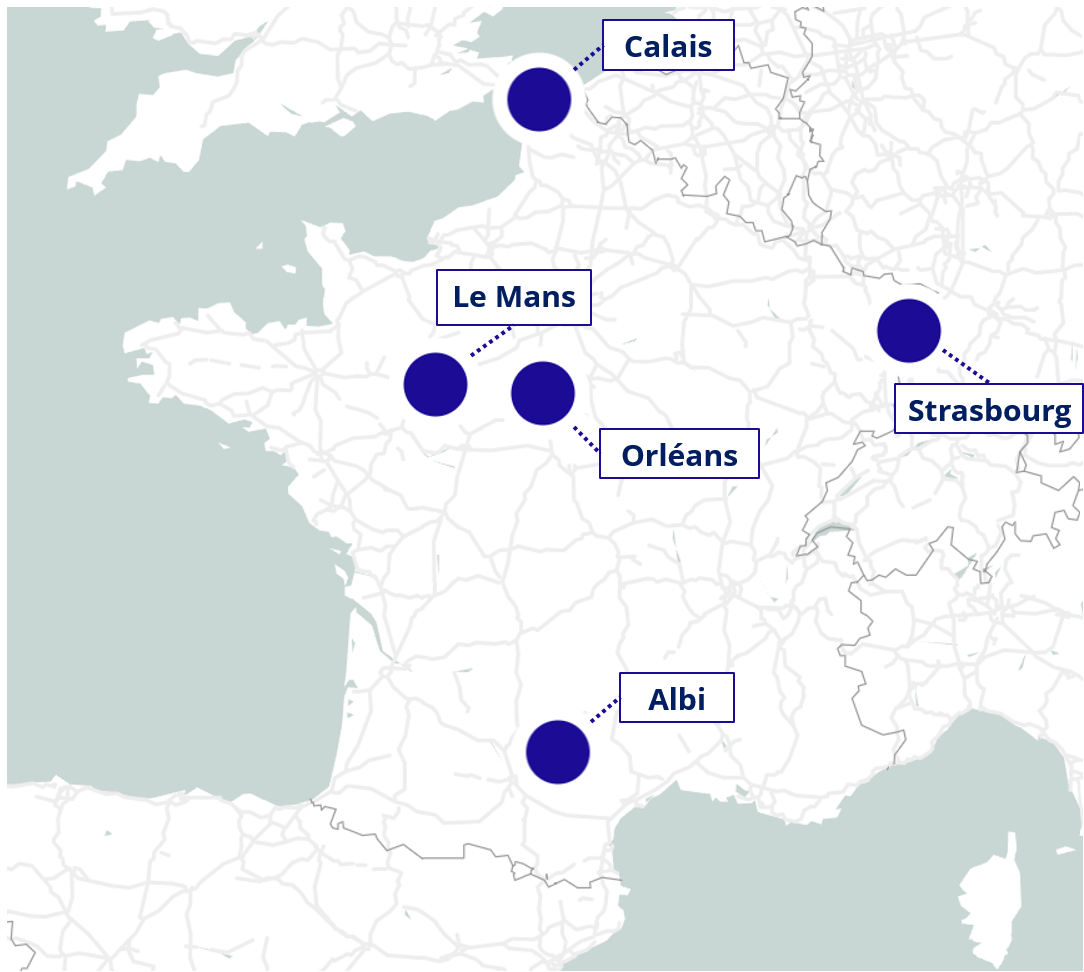}
    \caption{PPA locations.}
    \label{fig:map}
\end{wrapfigure}
In order to fix PPA prices for the model, first a worst-case (highest cost) LCOE for each site was calculated from the lowest yearly capacity factor (CF), as shown in Table \ref{tab:ppa_prices}. The LCOE's were used as a guide to establish a rough ordering of PPA prices, which were then chosen arbitrarily to provide a range of reasonably competitive PPA options at levels coherent with recent average PPA prices reported for the European market (\cite{Joinville2023European3Q}).

A summary of capacity factors and prices for the PPA's used is given in Table \ref{tab:ppa_prices}. A multivariate normal distribution is obtained from these historical capacity factors, and this is sampled randomly to obtain new capacity factors during scenario generation. 

\begin{table}
\centering
\scriptsize
\begin{tabular}{|c|c|p{0.08\linewidth}|p{0.08\linewidth}|p{0.08\linewidth}|p{0.08\linewidth}|p{0.09\linewidth}|}
\hline
\multicolumn{1}{|c|}{\textbf{PPA}} & \textbf{Location} & \textbf{CF \,Min} & \textbf{CF Mean} & \textbf{CF \,Max} & \textbf{Max LCOE} & \textbf{Chosen PPA Price}     \\ \hline
\multirow{4}{*}{\rotatebox[origin=c]{90}{\tiny \textbf{Solar PV}}} & Le Mans           & 14.5 \%         & 15.0 \%          & 16.2 \%         & 58                & 71   \\ \cline{2-7} 
                                   & Calais            & 14.2 \%         & 14.8 \%          & 15.4 \%         & 59                & 72 \\ \cline{2-7} 
                                   & Strasbourg        & 14.7 \%         & 15.3 \%          & 16.2 \%         & 57                & 68  \\ \cline{2-7} 
                                   & Albi              & 16.0 \%         & 16.8 \%          & 17.5 \%         & 52                & 66  \\ \hline
\multirow{5}{*}{\rotatebox[origin=c]{90}{\tiny \textbf{Wind}}}     & Orleans           & 25.5 \%         & 27.7 \%          & 31.4 \%         & 54                & 70 \\ \cline{2-7} 
                                   & Le Mans           & 27.1 \%         & 29.3 \%          & 32.7 \%         & 51                & 69  \\ \cline{2-7} 
                                   & Calais            & 40.9 \%         & 42.7 \%          & 46.5 \%         & 34                & 65    \\ \cline{2-7} 
                                   & Strasbourg        & 15.1 \%         & 16.9 \%          & 19.5 \%         & 92                & 80  \\ \cline{2-7} 
                                   & Albi              & 24.6 \%         & 25.8 \%          & 26.5 \%         & 56                & 78   
\\ \hline
\end{tabular}
\caption{PPA details: capacity factors (CF) for 2013-2022, maximum LCOE in €/MWh corresponding to the worst capacity factor, price chosen arbitrarily for the case study in €/MWh.}
\label{tab:ppa_prices}
\end{table}

\subsubsection{Green Hydrogen Subsidies}\label{subsubsec:case_study:components_uncertainty:green_subsidy}
Green hydrogen subsidies when applied in this paper are modelled using a requirement for hourly-time matching of hydrogen production with renewable sources. This rule is to be applied in the US and Europe post-2030 (\cite{Guillotin2025HydrogenComparison}).

A subsidy level of 3 €/kg is chosen arbitrarily for analysis but is considered reasonable, given the 3 \$/kg maximum tax incentive available in the U.S. and the 4.5 €/kg subsidy cap in the European Hydrogen Bank auctions (\cite{EuropeanCommission2023CommissionProduction.,USDepartmentofEnergy2023AssessingCredit}), although it should be noted that the first round of the European auctions cleared much lower, at 0.48 €/kg (\cite{EuropeanCommission2024EuropeanEurope}).

\subsubsection{Cost of Capital}\label{subsubsec:case_study:components_uncertainty:coc}
A fixed cost of capital of 5\% is used throughout. This is chosen as a reasonable value based on the contextual assumptions of this study --- that the project takes place in Europe and with a long-term HPA for 100\% of intended production (\cite{InternationalEnergyAgency2021TheTransitions}). If merchant risk were involved for hydrogen sales, a larger cost of capital would be appropriate.  

\subsubsection{Scenario Sets}\label{subsubsec:case_study:components_uncertainty:scenario_gen}
The deterministic Policies 1 and 2 of Section \ref{subsec:policies:defintions} were optimised on 1 expected value scenario. The stochastic policies 3, 4 and 5 were optimised using 25 randomly sampled scenarios (see the Supplementary Material for validation of this number of scenarios). 

The solutions of all policies are then tested on an out-of-sample set of 1000 scenarios. 

A strict separation between in-sample and out-of-sample scenarios is applied: even-numbered years are used to create the in-sample scenarios, and odd-numbered years are used to create the out-of-sample scenarios. Details of the scenario generation process are provided in the Supplementary Material.

\begin{table}
\centering
\scriptsize
\begin{tabular}{llcl}
\hline
\multicolumn{4}{c}{\textbf{Case Study Parameters}}                                                                                                                                                                                      \\ \hline
\multicolumn{1}{c}{\multirow{2}{*}{\begin{tabular}[c]{@{}c@{}}Electrolyser + \\ Compressor\end{tabular}}}   & Combined CAPEX                                                        & 1.7 M                     & €/MW (Elec)           \\
\multicolumn{1}{c}{}                                                                                        & Combined Efficiency                                                   & 56 \%                     & MWh (H2) / MWh (Elec) \\
                                                                                                            & Lifetime                                                              & 13                        & Years                 \\
                                                                                                            &                                                                       &                           &                       \\
\multicolumn{1}{c}{H2 Storage}                                                                              & CAPEX (Energy)                                                        & 75k                       & €/MWh                 \\
                                                                                                            & CAPEX (Power)                                                         & 50k                       & €/MW                  \\
                                                                                                            & Energy Efficiency                                                     & 100 \%                    &                       \\
                                                                                                            & Power Efficiency                                                      & 100 \%                    &                       \\
                                                                                                            & Usable Capacity                                                       & 100 \%                    &                       \\
                                                                                                            & Initial SOC                                                           & 50 \%                     &                       \\
                                                                                                            & Lifetime                                                              & 25                        & Years                 \\
                                                                                                            &                                                                       & \multicolumn{1}{l}{}      &                       \\
\multicolumn{1}{c}{Network Connection}                                                                      & CAPEX                                                                 & 75k                       & €/MW                  \\
                                                                                                            & Lifetime                                                              & 25                        & Years                 \\
                                                                                                            &                                                                       & \multicolumn{1}{l}{}      &                       \\
\multicolumn{1}{c}{\multirow{2}{*}{\begin{tabular}[c]{@{}c@{}}Project \\ (Global Parameters)\end{tabular}}} & Discount Rate                                                         & 5 \%                      &                       \\
\multicolumn{1}{c}{}                                                                                        & Lifetime                                                              & 25                        & Years                 \\
                                                                                                            & H2 Energy Density                                                     & 33.33                     & kg H2 / MWh H2        \\
                                                                                                            & \begin{tabular}[c]{@{}l@{}}Demand Curtailment \\ Penalty\end{tabular} & \multicolumn{1}{l}{1,000} & €/MWh                 \\ \hline
\end{tabular}
\caption{Case study -- common parameters. Sources: \cite{CommissiondeRegulationdelEnergie2019CoutsContinentale,InternationalEnergyAgency2021TheTransitions,Wu2022ElectricityMarkets,Matute2023Techno-economicPPAs,Jodry2023IndustrialOptimisation}.}
\label{tab:paramaters}
\end{table}

\begin{table}[]
\tiny
\begin{tabular}{lllll}
\hline
\multicolumn{5}{l}{\textbf{Notation}}                                                                                                                          \\ \hline
\multicolumn{2}{l}{\textit{Sets}}                            &           & \multicolumn{2}{l}{\textit{Hedging Parameters}}                                     \\
$\mathbb{Q}$           & Futures product periods             &           & $H_q$                  & Futures product $q$ duration                               \\
                       &                                     &           &                        &                                                            \\
$\mathbb{A}$           & Power Purchase Agreements           &           & $m^{q}_h$              & Futures activation function                                \\
                       &                                     &           &                        &                                                            \\
$S$                    & Number of scenarios                 &           & $p_h$                  & Futures peakload activation function                       \\
                       &                                     &           &                        &                                                            \\
$H$                    & Number of hours                     &           & \multicolumn{2}{l}{\textit{Technical Parameters}}                                   \\
                       &                                     &           & $\eta^{ez}$            & Electrolyser + Compressor energy                           \\
\textit{Indices}       &                                     &           &                        & efficiency (MWh H2 / MWh elec)                             \\
$s$                    & Scenario index                      &           &                        &                                                            \\
                       &                                     &           & $\eta^{hs-p,in}$       & Hydrogen storage input power                               \\
$h$                    & Hour of the year index              &           &                        & efficiency                                                 \\
                       &                                     &           &                        &                                                            \\
$i$                    & Network access time-slot index      &           & $\eta^{hs-p,out}$      & Hydrogen storage discharge power                           \\
                       &                                     &           &                        & efficiency                                                 \\
$q$                    & Futures delivery period index       &           &                        &                                                            \\
                       &                                     &           & \textit{$\eta^{hs-e}$} & Hydrogen storage energy                                    \\
$a$                    & PPA contract index                  &           &                        & efficiency                                                 \\
                       &                                     &           & \textit{}              &                                                            \\
\multicolumn{2}{l}{\textit{Design Decisions}}                &           & \multicolumn{2}{l}{\textit{Cost Parameters}}                                        \\
$x^{ez-p}$             & Electrolyser power (MW elec)        &           & $C^{ez}$               & Electrolyser CAPEX                                         \\
                       &                                     &           & \textit{}              & (€/MW elec/year)                                           \\
$x^{hs-e}$             & Hydrogen storage energy (MWh H2)    &           & \textit{}              &                                                            \\
                       &                                     &           & $C^{hs-e}$             & Hydrogen Storage Energy                                    \\
$x^{hs-p}$             & Hydrogen storage power (MW H2)      &           & \textit{}              & CAPEX (€/MWh H2/year)                                      \\
                       &                                     &           &                        &                                                            \\
$x^{nw-p}$             & Network connection power (MW)       &           & $C^{hs-p}$             & Hydrogen Storage Power                                     \\
                       &                                     &           &                        & CAPEX (€/MW H2/year)                                       \\
$x^{BaF-e,q}$          & Baseload futures, purchased energy  &           & \textit{}              &                                                            \\
                       & for period $q$ (MWh)                &           & \textit{$C^{nw-p}$}    & Network Connection Power                                   \\
                       &                                     &           &                        & CAPEX (€/MW/year)                                          \\
$x^{PkF-e,q}$          & Peakload futures, purchased energy  &           &                        &                                                            \\
                       & for period $q$ (MWh)                &           & $C^{BaF-e,q}$          & Baseload futures price,                                    \\
                       &                                     &           &                        & period $q$ (€/MWh)                                         \\
$\hat{x}^{nw,i}$       & Network access power subscription   &           &                        &                                                            \\
                       & for time-slot $i$ (MW)              &           & $C^{PkF-e,q}$          & Peakload futures price,                                    \\
                       &                                     &           &                        & period $q$ (€/MWh)                                         \\
$x^{ppa-p,a}$          & PPA peak power for park $a$ (MWp)   &           &                        &                                                            \\
                       &                                     &           & $B^{nw,i}$             & Network power subscription cost                            \\
\multicolumn{2}{l}{\textit{Operational Decisions}}           &           &                        & for time-slot $i$ (€/MW/year)                              \\
$u^{da,out}_{s,h}$     & Spot market energy bought (MWh)     &           &                        &                                                            \\
                       &                                     &           & $G^{ppa-e,a}$          & PPA energy cost for park $a$ (€/MWh)                       \\
$u^{da,in}_{s,h}$      & Spot market energy sold (MWh)       & \textit{} &                        &                                                            \\
                       &                                     &           & $G^{hd-curt}$          & Hydrogen demand curtailment                                \\
$u^{ez,in}_{s,h}$      & Electrolyser dispatching (MW elec)  & \textit{} &                        & penalty (€/MWh H2)                                         \\
                       &                                     &           &                        &                                                            \\
$u^{hs,in}_{s,h}$      & Hydrogen storage input (MW H2)      &           & $G^{r-e}$              & Green hydrogen subsidy (€/MWh H2)                          \\
                       &                                     &           &                        &                                                            \\
$u^{hs,out}_{s,h}$     & Hydrogen storage output (MW H2)     &           & \multicolumn{2}{l}{\textit{Design Upper Bounds (UB)}}                               \\
                       &                                     &           & $X^{ez-p}$             & Electrolyser power UB (MW elec)                            \\
$u^{ppa,out,a}_{s,h}$  & PPA dispatched power, park $a$ (MW) &           &                        &                                                            \\
                       &                                     &           & $X^{hs-e}$             & Hydrogen storage energy                                    \\
$u^{ppa,curt,a}_{s,h}$ & PPA curtailed power, park $a$ (MW)  & \textit{} &                        & UB (MWh H2)                                                \\
                       &                                     &           & \textit{}              & \textit{}                                                  \\
$u^{hd,in}_{s,h}$      & Hydrogen demand fulfilled (MW H2)   & \textit{} & $X^{hs-p}$             & Hydrogen storage power UB (MW H2)                          \\
                       &                                     &           &                        &                                                            \\
$u^{hd,curt}_{s,h}$    & Hydrogen demand unmet (curtailed)   & \textit{} & $X^{nw-p}$             & Network connection power UB (MW)                           \\
                       & (MW H2)                             &           &                        &                                                            \\
                       &                                     &           & $X^{BaF-e,q}$          & Baseload futures, purchased energy                         \\
$u^{nw}_{s,h}$         & Power transferred to the site       &           &                        & for period $q$ UB (MWh)                                    \\
                       & from the network (MW)               &           &                        &                                                            \\
                       &                                     &           & $X^{PkF-e,q}$          & Peakload futures, purchased energy                         \\
$\delta_{s,h}$         & Subsidy-eligible hydrogen           &           &                        & for period $q$ UB (MWh)                                    \\
\textit{}              & produced (MWh H2)                   &           &                        &                                                            \\
                       &                                     &           & $\hat{X}^{nw,i}$       & Network access power subscription                          \\
\multicolumn{2}{l}{\textit{Operational States}}              &           &                        & for time-slot $i$ UB (MW)                                  \\
$z^{hs}_{s,h}$         & Hydrogen storage state-of-charge    &           &                        &                                                            \\
                       & (MWh)                               &           & $X^{ppa-p,a}$          & PPA peak power contracted                                  \\
\textit{}              &                                     &           &                        & for park $a$ UB (MWp)                                      \\
\multicolumn{2}{l}{\textit{Operational Uncertainties}}       &           &                        &                                                            \\
$w^{da}_{s,h}$         & Day-ahead market price (€/MWh)      &           & \multicolumn{2}{l}{\textit{Objective Function Paramaters}}                          \\
                       &                                     &           & $\alpha$               & CVaR quantile setting \{\{{[}\}\}0:1)                      \\
$w^{ppa,a}_{s,h}$      & PPA availability factor             &           &                        &                                                            \\
                       & (\% of peak power)                  &           & \textit{$\beta$}       & \textit{Risk aversion weighting \{\{{[}\}\}0:1\{\{{]}\}\}} \\
                       &                                     &           & \textit{}              &                                                            \\
$w^{hd,in}_{s,h}$      & Hydrogen offtaker demand (MW H2)    &           & \textit{}              &                                                            \\ \hline
\end{tabular}
\caption{Mathematical notation used. The designations $p$ and $e$ denote power and energy values respectively. UB = Upper bound.}
\label{tab:notation}
\end{table}

\subsection{Mathematical Model} \label{subsec:case_study:model}
The following equations form a linear program of the model represented in Figure \ref{fig:schema}. The variables used are summarised in Table \ref{tab:notation}.

Their feasible space of design decisions $\mathbb{X}$ is defined for the case study by the bound constraint Equations (\ref{eq:design_bounds_first})-(\ref{eq:design_bounds_last}). For each scenario $s$, operational (second stage) uncertainties are revealed (day-ahead price $w^{da}_{s,h}$, PPA production $w^{ppa,a}_{s,h}$ and hydrogen demand $w^{hd,in}_{s,h}$). The hourly operation (dispatching) of assets is solved deterministically within each scenario, in the feasible space $\mathbb{U}_{s,h}$ which is defined for the case study by Equations (\ref{eq:op_bounds_first})-(\ref{eq:op_bounds_last}). The only operational state variable used is hydrogen storage state-of-charge, for which the feasible space $\mathbb{Z}_{s,h}$ is constrained by Equations (\ref{eq:states_first})-(\ref{eq:states_last}). Design cost $J^d$ is defined by Equation (\ref{eq:design_costs}), and operational costs $J^{o}_{s}$ are defined by Equations (\ref{eq:op_cost_start})-(\ref{eq:op_cost_end}).

\subsubsection{Design Decision Bound Constraints}\label{subsec:model:bounds}
The following constraints impose the design decision upper bounds contained in the vector $X$, which may represent limited availability of contracts, limited space to install physical assets, limited network capacity, etc:

\begin{align} 
0 &\leq x^{ez-p} \leq X^{ez-p},  \label{eq:design_bounds_first} \\
0 &\leq x^{hs-e} \leq X^{hs-e},  \\
0 &\leq x^{hs-p} \leq X^{hs-p},  \\
0 &\leq x^{nw-p} \leq X^{nw-p},  \\
0 &\leq x^{BaF-e,q} \leq X^{BaF-e,q},  \\
0 &\leq x^{PkF-e,q} \leq X^{PkF-e,q},  \\
& \forall \quad q \in \mathbb{Q}, \nonumber \\  \noalign{\vskip12pt} 
\text{and } \quad 0 &\leq x^{ppa-p,a} \leq X^{ppa-p,a} , \label{eq:design_bounds_last}\\
& \forall \quad a \in \mathbb{A}. \nonumber
\end{align}

Upper bounds are never reached for any design decision in this paper.

\subsubsection{Operational Decision Limit Constraints}\label{subsec:model:limits}
The operational decisions $u$ are limited by the following constraints for all $h \in [1..H]$, and all $s \in [1..S]$:

\begin{align} 
0 &\leq u^{ez,in}_{s,h} \leq {x}^{ez-p}, \label{eq:op_bounds_first}\\
0 &\leq u^{hs,in}_{s,h} \leq {x}^{hs-p}, \\
0 &\leq u^{hs,out}_{s,h} \leq {x}^{hs-p}, \\
0 &\leq u^{nw}_{s,h} \leq {x}^{nw-p}, \\
0 &\leq u^{hd,in}_{s,h} \leq {w}^{hd}_{s,h}, \\
0 &\leq u^{hd,curt}_{s,h} \leq {w}^{hd}_{s,h}, \\
 \noalign{\vskip12pt}
0 &\leq u^{ppa,a,out}_{s,h} \leq {w}^{ppa,a}_{s,h} \cdot x^{ppa-p,a}, \\
\text{and } \quad 0 &\leq u^{ppa,a,curt}_{s,h}  \leq {w}^{ppa,a}_{s,h} \cdot x^{ppa-p,a}, \\
& \forall \quad a \in \mathbb{A}. \nonumber 
\end{align}

\subsubsection{Energy Balance Constraints} \label{subsubsec:case_study:model:constraints}
The following operational constraints ensure an energy balance at each node for all $h \in [1..H]$ and all $s$ in $[1..S]$:\\

\small \noindent \textit{Hydrogen Balance Constraint}

\noindent The electrolyser is modelled as a single linear efficiency $\eta^{ez}$. The energy balance at the hydrogen node is given as: 
\begin{equation} \label{eq:h2_balance}
u^{ez,in}_{s,h} \cdot \eta^{ez} + u^{hs,out}_{s,h} = u^{hs,in}_{s,h} +  w^{hd,in}_{s,h} - u^{hd,curt}_{s,h}.
\end{equation}

\small \noindent \textit{Electricity Balance Constraint}

Next, electricity delivered to the electrolyser is transferred across the connection to the larger network:
\begin{align}\label{eq:elec_balance}
& u^{nw}_{s,h} = u^{ez,in}_{s,h}.
\end{align}

\small \noindent \textit{Market Balance Constraint}

\noindent Within the larger network, a balanced energy market perimeter is assumed (i.e., energy bought is balanced with energy consumed and sold): 
\begin{align} \label{eq:market_balance}
u^{da,out}_{s,h} + \sum_{a=1}^{A} u^{ppa,out,a}_{s,h} \quad \quad \quad \quad \quad \quad \quad \quad \quad \quad \quad  &  \\
\quad + \sum_{q \in \mathbb{Q}} m^{q}_{h} \cdot \{x^{BaF-e,q}/H_q + p_h \cdot x^{PkF-e,q}/(H_q/2)\} & \nonumber \\
& = u^{da,in}_{s,h} + u^{nw}_{s,h}, \nonumber \\ \noalign{\vskip6pt}
\text{ where } \quad p_h = \pmb1_{\{8 \leq h \: \text{mod} \: 24 < 20\}} \quad \quad \quad \quad \quad \quad &  \\ \noalign{\vskip6pt}
\text{ and } m^{q}_h = {\pmb1_{\{H^q_{begin} \leq h \leq H^q_{end}\}}}. \quad \quad \quad \quad \quad \quad & 
\end{align}
The unit function $p_h$ activates Peakload Futures purchases only during the hours 8am to 8pm, and similarly, the unit function $m^{q}_h$ ensures that futures product $q$ is only activated during its delivery period (futures prices and parameters are provided in the Supplementary Material).\\

\small \noindent \textit{PPA Internal Balance Constraint}

\noindent PPA production can also be curtailed (it is assumed that the respective park is able to curtail its production in the case of negative day-ahead prices, for example):
\begin{align} \label{eq:ppa_intern_balance}
& u^{ppa,out,a}_{s,h} = w^{ppa,a}_{s,h} \cdot x^{ppa-p,a} - u^{ppa,curt,a}_{s,h} \\
& \forall \quad a \in \mathbb{A}. \nonumber
\end{align}
Note that any curtailed energy is still ‘paid-for' in the operational cost, thus the formulation is akin to a \textit{take-or-pay} type PPA.\\

\subsubsection{No-Resale (NR) Constraint}\label{subsec:model:no_arbitrage}

As discussed in Section \ref{subsubsec:intro:lit:noresale}, some planning policies use a no-resale (NR) formulation. In this case, an additional constraint is added which eliminates resale of excess energy on the day-ahead market:
\begin{equation} \label{eq:op_bounds_last}
u^{da,in}_{s,h} = 0.
\end{equation}
Energy resale is always allowed for testing (i.e., Equation (\ref{eq:op_bounds_last}) does not apply in out-of-sample testing for any planning policy).

\subsubsection{Storage State-of-Charge}\label{subsec:model:storage}

The hydrogen storage unit state-of-charge (SOC) is modelled with the following state constraint:
\begin{align} \label{eq:states_first}
& z^{hs}_{s,h+1} = z^{hs}_{s,h} \cdot (1-\eta^{hs-e} )  + (\eta^{hs-p,in} \cdot u^{hs,in}_{s,h}  - \frac{ u^{hs,out}_{s,h}} {\eta^{hs-p,out}}).
\end{align}
The charging and discharging efficiencies are given by $\eta^{hs-p,in}$ and $\eta^{hs-p,out}$, and $\eta^{hs-e}$ is the energy storage auto-discharge loss rate per hour. The maximum SOC is limited by the design decision for storage rated capacity:

\begin{align}
0 &\leq z^{hs}_{s,h} \leq {x}^{hs-e}.
\end{align}

\noindent The initial state of charge is assumed to be half of rated capacity: 
\begin{align} \label{eq:soc_init}
& z^{hs}_{s,0} = x^{hs-e} \cdot 0.5.
\end{align}

\noindent And the final state of charge at the end of the year is required to be equal to or greater than the initial capacity:
\begin{align} \label{eq:states_last}
& z^{hs}_{s,H} \geq z^{hs}_{s,0}.
\end{align}

\subsubsection{Design Costs}\label{subsec:model:design_decisions}
The design decisions attract design costs based on their respective price parameters in the vector $C$:

\begin{align} \label{eq:design_costs}
J^d(x,C,B) \quad =& \quad \quad x^{ez-p} \cdot C^{ez-p}  + x^{hs-e} \cdot C^{hs-e} \\
& \quad +  x^{hs-p} \cdot C^{hs-p} + x^{nw-p} \cdot C^{nw-p} \nonumber\\
& \quad + \sum_{q \in \mathbb{Q}} \{x^{BaF-e,q} \cdot C^{BaF-e,q} + x^{PkF-e,q} \cdot C^{PkF-e,q}\} \nonumber\\
& \quad + J^{d,nw}(x,B). \nonumber
\end{align}
Network access peak power charges $J^{d,nw}(x,B)$ are defined as a linear function of the design decisions and network access cost parameters $B$ (additional details of the network access charge regime used are provided in the Supplementary Material). \\

CAPEX costs are annualised using the following formula, where $LT^k$ is the lifetime of asset $k$ and $\tau$ is the discount rate: 
\begin{align} \label{eq:discount_capex}
    & C^k = CAPEX^k \cdot \frac{\tau \cdot (\tau + 1)^{LT^k}}{(\tau + 1)^{LT^k} - 1}.
\end{align}

\subsubsection{Operational Costs}\label{subsec:model:op}
Operational Costs are the net sum of revenues and costs on the day-ahead market, power purchased from PPA contracts, penalties paid for non-served demand, network access energy costs, and green hydrogen subsidies. 

Fixed parameters relating to operational costs include the energy cost for each PPA $G^{ppa-e,a}$, the cost of penalties for curtailed hydrogen demand with the flat cost $G^{hd-curt}$, network access time-of-use energy charges with the cost vector $G^{nw-e}_h$, and the flat-rate green hydrogen subsidy $G^{r-e}$:
\begin{align} \label{eq:op_cost_start}
J^o_s(x,u_s,w_s,G) = \sum_{h=1}^{H}& 
\left\{
    \begin{array}{lr}
        (u^{da,out}_{s,h} - u^{da,in}_{s,h}) \cdot w^{da}_{s,h} \\ \noalign{\vskip6pt} 
        + \sum_{a \in \mathbb{A}} \{x^{ppa-p,a} \cdot w^{ppa-p,a}_{s,h} \cdot G^{ppa-e,a}\} \\ \noalign{\vskip6pt}
       + u^{hd,curt}_{s,h} \cdot G^{hd-curt} + u^{nw}_{s,h} \cdot G^{nw-e}_h -
        \delta_{s,h} \cdot G^{r-e}
    \end{array}
\right\} \\ \noalign{\vskip6pt} \nonumber
& \forall \quad s \in [1..S]. \\ 
 \nonumber
\end{align}

Where hydrogen subsidies are not included, $G^{r-e}=0$ in Equation (\ref{eq:op_cost_start}).

\subsubsection{Green Hydrogen Subsidies}\label{subsec:model:rfnbo}

Green hydrogen subsidies are assumed to be only attainable with the use of hour-by-hour temporally correlated renewable power (in the context of this study, this means power originating from renewable PPA's). The subsidy-eligible quantity of hydrogen $\delta_{s,h}$ produced at hour $h$ in scenario $s$, is defined with the constraints:

\begin{align} \label{eq:rfnbo_subsidy}
& \delta_{s,h}  \leq  (\sum_{a \in \mathbb{A}} u^{ppa,out,a}_{s,h} ) \cdot \eta^{ez},  \\
& \delta_{s,h}  \leq  u^{ez,in}_{s,h} \cdot \eta^{ez},  \label{eq:op_cost_end}\\
& \forall \quad h \in [1..H], \quad s \in [1..S].  \nonumber
\end{align}

These constraints effectively perform a \textit{max} operation between the total PPA power available in a given hour, and the electrolyser production in the same hour. Thus making sure that no more hydrogen can be counted than is actually produced for that hour, and that the amount of \textit{green} hydrogen counted in that hour cannot exceed the (conversion efficiency-adjusted) renewable power received for that hour from the PPA's. Additional power from PPA's which is instead sold to the market in a given hour thus forgoes its `green' credentials.

\section{Results and Discussion} \label{sec:results_methodo}
\begin{figure}
  \centering
\includegraphics[scale=0.43]{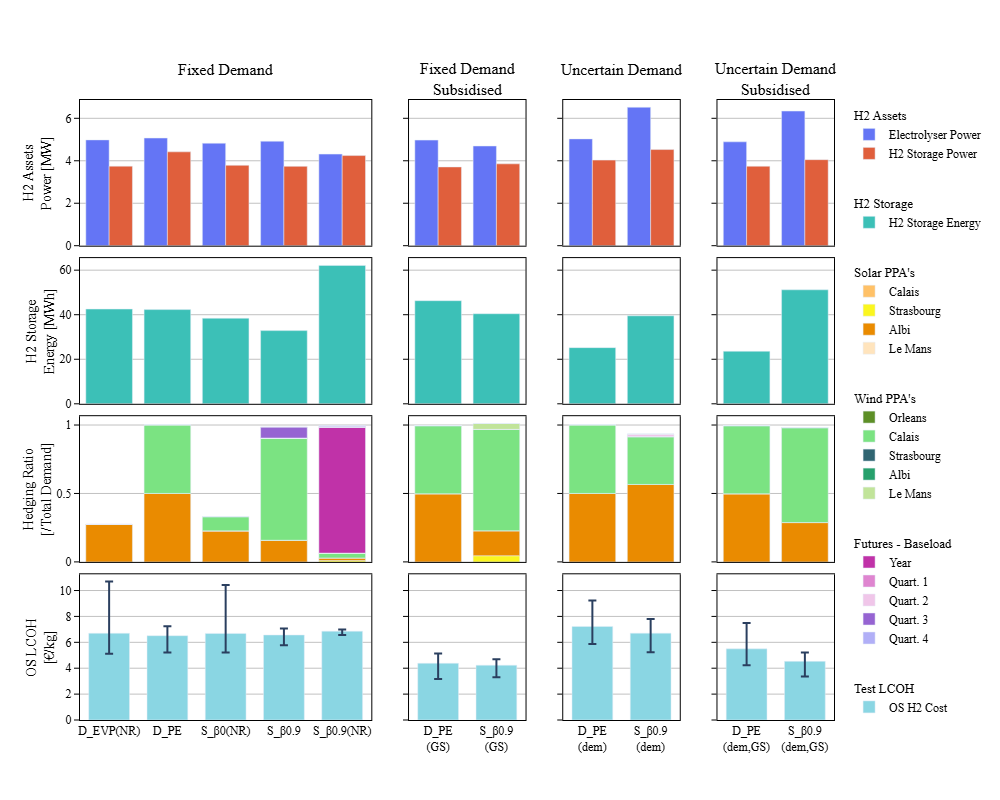}
    \caption{(Top row) Electrolyser and hydrogen storage input capacity (power); (Second row) H2 storage installed energy capacity; (Third row) hedging ratio of futures and Power Purchase Agreement (PPA) hedging decisions (hedging ratio is the proportion of the expected volume for the in-sample scenarios with respect to total yearly electricity demand); (Bottom row) Out-of-sample LCOH results. Bars indicate mean values across all scenarios. Whiskers indicate maximum and minimum scenario Levelised Cost of Hydrogen (LCOH).}
    \label{fig:results_plot}
\end{figure}

\begin{figure}
  \centering
\includegraphics[scale=0.35]{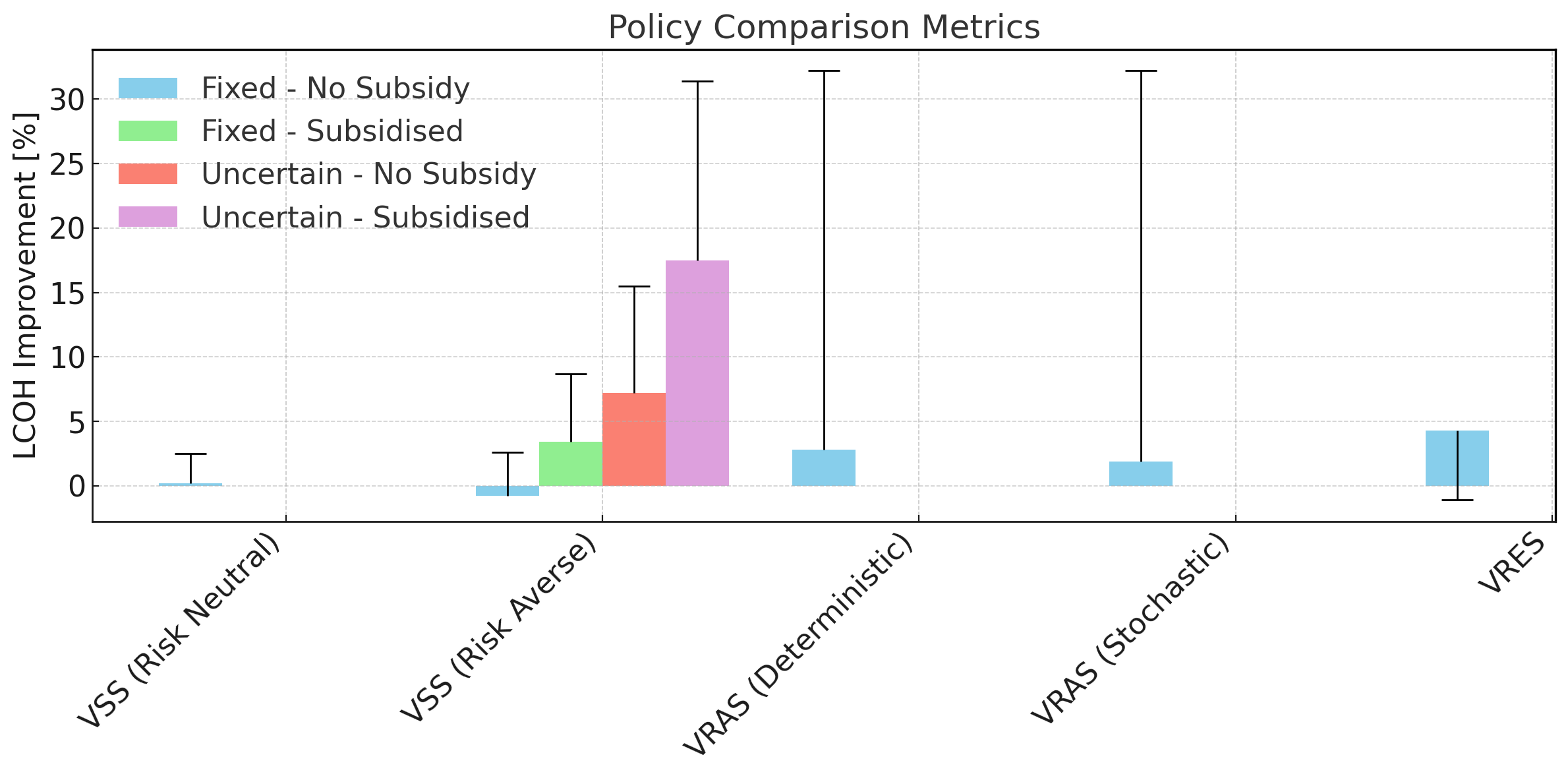}
    \caption{Policy Comparison Metrics results, also described in Table \ref{tab:policy_comparison}. Bars represent mean LCOH reduction. Whiskers represent the reduction in worst-case LCOH.}
    \label{fig:policy_comp_plot_1}
\end{figure}

The results are summarised in Figure \ref{fig:results_plot}, and the out-of-sample test results for each solution in Table \ref{tab:results_lcoh}. The computed Policy Comparison Metrics are shown in Table \ref{tab:policy_comparison} and displayed in Figures \ref{fig:policy_comp_plot_1} and \ref{fig:policy_comp_plot_2}. Additional data on equipment sizes and hedging decisions are provided in the Supplementary Material.

\begin{table}[]
\scriptsize
\begin{tabular}{@{}lcccccccc@{}}
\toprule
                              & \multicolumn{4}{c}{\textbf{Fixed Demand HPA}}                                     & \multicolumn{4}{c}{\textbf{Uncertain Demand HPA}}                                 \\ \midrule
\textbf{Metric}               & \multicolumn{2}{c}{\textbf{No Subsidy}} & \multicolumn{2}{c}{\textbf{Subsidised}} & \multicolumn{2}{c}{\textbf{No Subsidy}} & \multicolumn{2}{c}{\textbf{Subsidised}} \\ \midrule
\textit{VSS (Risk Neutral)}   & \textit{(a)}        & 0.2 / 2.5         &                      &                  &                     &                   &                     &                   \\
\textit{VSS (Risk Averse)}    & \textit{(b)}        & -0.8 / 2.6        & \textit{(f)}         & 3.4 / 8.7        & \textit{(g)}        & 7.2 / 15.5        & \textit{(j)}        & \textbf{17.5} / \textbf{31.4}       \\
\textit{VRAS (Deterministic)} & \textit{(c)}        & 2.8 / \textbf{32.2}        &                      &                  &                     &                   &                     &                   \\
\textit{VRAS (Stochastic)}    & \textit{(d)}        & 1.9 / \textbf{32.2}        &                      &                  &                     &                   &                     &                   \\
\textit{VRES}                 & \textit{(e)}        & 4.3 / -1.1        &                      &                  &                     &                   &                     &                   \\
\textit{CDU (Deterministic)}  &                     &                   &                      &                  & \textit{(h)}        & 11.0 / 27.4       & \textit{(k)}        & 20.2 / 31.4       \\
\textit{CDU (Stochastic)}     &                     &                   &                      &                  & \textit{(i)}        & 2.2 / 10.5        & \textit{(l)}        & 6.6 / 10.0        \\ \bottomrule
\end{tabular}
\caption{Policy Comparison Metrics: Value of the Stochastic Solution (VSS), Value of the Risk Averse Solution (VRAS), Value of the Re-sale Enabled Solution (VRES), Cost of Demand Uncertainty (CDU). Values are provided in percentage (\%) difference between solution LCOH performance and are displayed in the format {[} Mean / Max {]}. Specific solutions compared for each result are: \textit{(a) D\_EVP(NR) -- S\_$\beta$0(NR); (b) D\_PE -- S\_$\beta$0.9; (c) D\_EVP(NR) -- D\_PE; (d) S\_$\beta$0(NR) -- $S\_\beta0.9$; (e) S\_$\beta$0.9(NR) -- $S\_\beta0.9$; (f) D\_PE(dem) -- S\_$\beta$0.9(dem); (g) $D\_PE(dem)$} --~\textit{$S\_\beta0.9(dem)$; (h) $D\_PE(dem)$} -- ~\textit{$D\_PE$; (i) $S\_\beta0.9(dem)$} --~\textit{$S\_\beta0.9$; (j) $D\_PE(dem,GS)$} --~\textit{$S\_\beta0.9(dem,GS)$; (k) $D\_PE(dem,GS)$} -- ~\textit{$D\_PE(GS)$; (l) $S\_\beta0.9(dem,GS)$} --~\textit{$S\_\beta0.9(GS)$}.}
\label{tab:policy_comparison}
\end{table}

\begin{wraptable}{r}{7cm}
\scriptsize
\centering
\begin{tabular}{lcc}
\hline
\textbf{}             & \multicolumn{2}{l}{\textbf{LCOH (€/kg)}}                                              \\ \cline{2-3} 
\textbf{Solution}     & \textbf{Mean}        & \textbf{\begin{tabular}[c]{@{}c@{}}Worst \\ Case\end{tabular}} \\ \hline
$D\_EVP(NR)$          & 6.70                 & 10.69                                                          \\
$D\_PE$               & 6.51                 & 7.25                                                           \\
$S\_\beta0(NR)$       & 6.69                 & 10.42                                                          \\
$S\_\beta0.9$         & 6.57                 & 7.06                                                           \\
$S\_\beta0.9(NR)$     & 6.86                 & 6.99                                                           \\
                      & \multicolumn{1}{l}{} & \multicolumn{1}{l}{}                                           \\
$D\_PE(GS)$           & 4.38                 & 5.14                                                           \\
$S\_\beta0.9(GS)$     & 4.23                 & 4.70                                                           \\
                      & \multicolumn{1}{l}{} & \multicolumn{1}{l}{}                                           \\
$D\_PE(dem)$          & 7.23                 & 9.24                                                           \\
$S\_\beta0.9(dem)$    & 6.71                 & 7.80                                                           \\
                      & \multicolumn{1}{l}{} & \multicolumn{1}{l}{}                                           \\
$D\_PE(dem,GS)$       & 5.49                 & 7.50                                                           \\
$S\_\beta0.9(dem,GS)$ & 4.53                 & 5.22                                                           \\ \hline
\end{tabular}
\caption{Out-of-sample LCOH results.}
\label{tab:results_lcoh}
\end{wraptable}

\subsection{Fixed Demand HPA} \label{subsec:results_methodo:fixed_demand}
As was the case for \citet{Palmer2023Risk-conscious.}, the low VSS indicates little to no value is obtained using a stochastic model for the Fixed Demand HPA. 

On the other hand, with the very significant VRAS, the advantage of the risk-averse policies with respect to the risk-neutral policies is clear. The Stochastic Risk Averse policy $S\_\beta0.9$ improves the mean performance of the risk-neutral $S\_\beta0$ by 1.9\% on average, and by 32\% for the worst-case scenario performance. The rule-based risk-averse Pessimistic Expert policy $D\_PE$ has a similar performance improvement over the deterministic risk-neutral Expected Value Problem policy. The solutions for the risk-averse policies $D\_PE$, $S\_\beta0.9$, and $S\_\beta0.9(NR)$ look to avoid high day-ahead market price exposure, and so opt for significantly greater hedging ratios than the risk-neutral policies' solutions $D\_EVP(NR)$ and $S\_\beta0(NR)$.

The Value of the Resale-Enabled Solution (VRES) is positive on mean performance (4.3\%) and slightly negative on worst-case performance (-1.1\%). This is to be expected, because while both policies $S\_\beta0.9$ and $S\_\beta0.9(NR)$ are equally risk-averse,  the no-resale formulation will systematically underestimate hedging options with good re-sale value, and so it finds less optimal solutions on average. Note that we are only able to obtain a VRES for the risk-averse case, as for the risk-neutral case energy resale cannot be enabled without rendering the results speculative and nonsensical, as discussed in Section \ref{sec:intro:speculation}. 

The Risk-Averse No-Resale policy $S\_\beta0.9(NR)$ makes markedly different hedging decisions to the other policies, opting for almost exclusively baseload futures rather than PPA's. We see that this change in procurement strategy is also linked to changes to hydrogen equipment decisions --- the flat delivery profile of power futures facilitates a smaller electrolyser sizing, thus increasing its load factor to 94\% (as opposed to 81\% for the other policies). This solution does, however, select a larger storage capacity in order to use hydrogen produced during the night to ride through daily peak demand. This illustrates that different procurement strategies necessitate different equipment sizing decisions, and the need for co-optimisation of these decisions.

\subsection{Fixed Demand HPA -- Subsidised} \label{subsec:results_methodo:case_2}
The added complexity of time-matched subsidies leads to some benefits in stochastic modelling. The VSS shows a 3\% reduction in average LCOH and a 9\% reduction in worst-case LCOH. 

The addition of the 3 €/kg green hydrogen subsidy achieves on average a 2.13 €/kg reduction in overall LCOH for the Pessimistic Expert policy $D\_PE(GS)$, and a 2.33 €/kg reduction for the Stochastic Risk Averse policy $S\_\beta0.9(GS)$. This difference is due to the fact that the Stochastic Risk Averse solution is able to produce green eligible hydrogen at 80\% of total yearly production compared to 75\% for the Pessimistic Expert solution.

The addition of green subsidies results in a significant increase in storage capacity --- especially for the Stochastic Risk-Averse solution $S\_\beta0.9(GS)$ which sees a 34\% increase relative to $S\_\beta0.9$. In regards to energy sourcing, the need for hourly time-matching leads to a need for complementary renewable production uncertainties: $S\_\beta0.9(GS)$ opts for two more expensive PPA's (wind at Le Mans for 69 €/MWh, and solar at Strasbourg for 68 €/MWh) in order to increase its likelihood of access to subsidies.

\subsection{Uncertain Demand HPA} \label{subsec:results_methodo:part_2}
With the presence of demand uncertainty, the VSS shows a clear improvement of 7\% on average and 16\% for the worst-case when using the Stochastic Risk-Averse policy. 

The value of stochastic modelling is further highlighted when comparing the Cost of Demand Uncertainty (CDU). The expected LCOH of the deterministic Pessimistic Expert solution when demand uncertainty is added increases by  11\% on average and 27\% for the worst-case. For the Stochastic Risk-Averse solution, the CDU is much smaller --- increasing the average LCOH by only 2.2\%, although the worst-case performance is still 10.5\% more expensive. 

The success of the Stochastic Risk-Averse solution in the presence of demand uncertainty is mostly thanks to increased asset sizing. A 33\% increase in electrolyser size is observed, going from 4.9 MW for the solution $S\_\beta0.9$ with Fixed Demand to 6.5 MW for $S\_\beta0.9(dem)$ with Uncertain Demand. This represents a reduction in electrolyser load factor to 61\%. Hydrogen storage energy capacity also increases by 20\%. 

While this increased asset sizing is necessary to accommodate longer periods of high demand, it does require a greater investment cost, with CAPEX increasing by 30\% for the solution $S\_\beta0.9(dem)$ with respect to $S\_\beta0.9$ of the Fixed Demand HPA. This higher CAPEX cost may be prohibitive for some project developers, as the additional financing may be difficult to obtain.
\begin{wrapfigure}{r}{5cm}
  \centering
\includegraphics[scale=0.35]{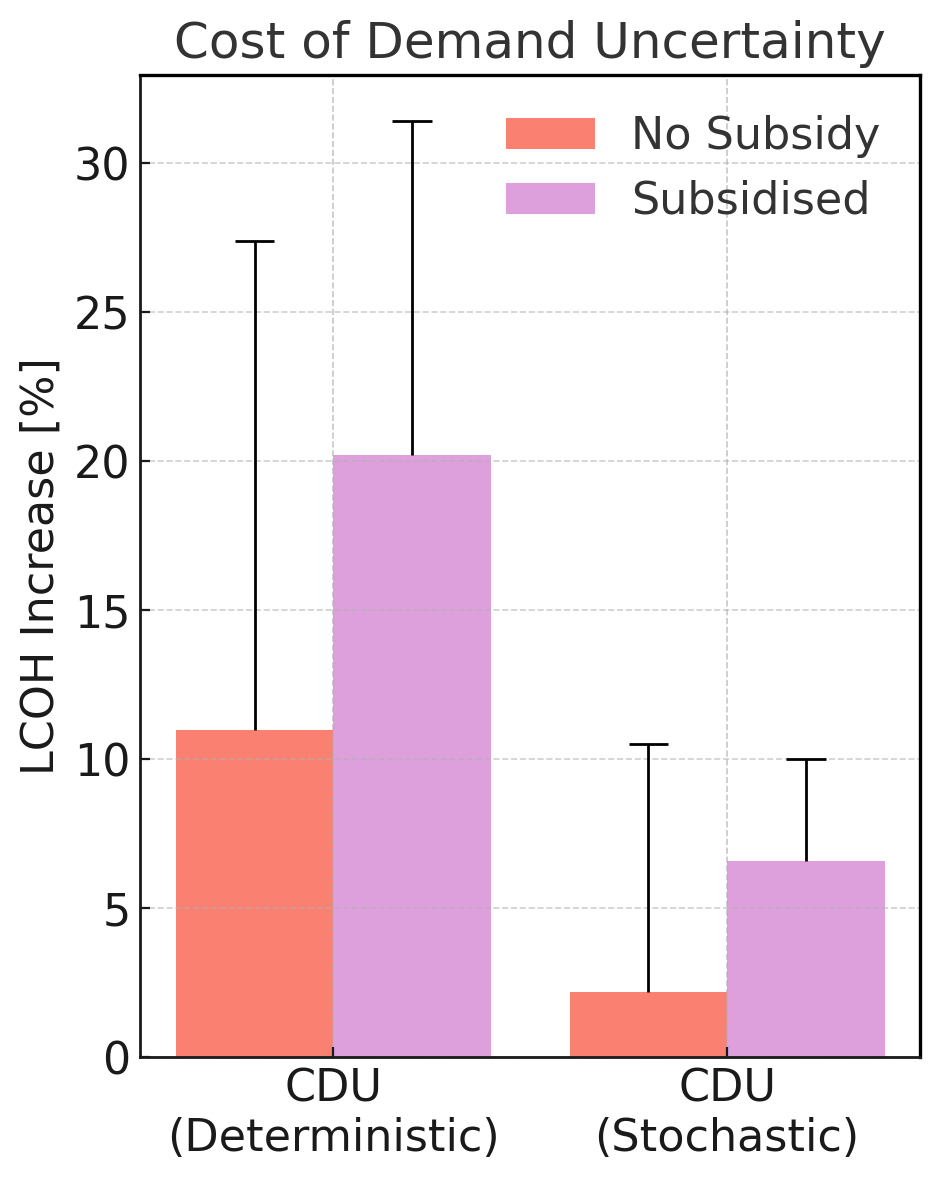}
    \caption{Cost of Demand Uncertainty results, also described in Table \ref{tab:policy_comparison}. Bars represent mean LCOH increase with demand uncertainty added. Whiskers represent the increase in worst-case scenario LCOH when demand uncertainty is added.}
    \label{fig:policy_comp_plot_2}
\end{wrapfigure}
\subsection{Uncertain Demand HPA -- Subsidised} \label{subsec:results_methodo:part_4}
In the context of both demand uncertainty and the possibility of green subsidies, the superiority of stochastic modelling is the strongest yet. The VSS in this case is such that the highest LCOH of 5.22 €/kg for the Stochastic Risk-Averse policy $S\_\beta0.9(dem,GS)$ is less than the mean LCOH of 5.49 €/kg for the Pessimistic Expert policy $D\_PE(dem,GS)$.

We also see that the Cost of Demand Uncertainty is increased when green subsidies are available: the cost of an Uncertain Demand HPA -- Subsidised is a 6.6\% increase on average with respect to a Fixed Demand HPA -- Subsidised (see Figure \ref{fig:policy_comp_plot_2}). 

For the Stochastic Risk-Averse policy, the possibility of green subsidies provokes a particularly large increase in hydrogen storage capacity, which increases by 29.6\% to 51.2 MWh for $S\_\beta0.9(dem,GS)$ with respect to 39.5 MWh for $S\_\beta0.9(dem)$. 

We also observe that for the Stochastic Risk-Averse policy, the addition of a green subsidy does not lead to the same increase in diversity of energy hedging as for the Fixed Demand HPA -- Subsidised. This is likely because achieving the necessary hourly time-matching by using increased asset sizing is more cost-effective than increasing hedging diversity, given that increased asset sizes are already necessary for managing the demand uncertainty.

\section{Conclusion and Policy Recommendations} \label{sec:conclusion}
Development of the low-carbon hydrogen economy is being slowed down by the difficulty in passing the final investment decision stage of the project cycle. This is in part due to the financial risk stemming from the sensitivity of production costs to a complex trinity of uncertainties in hydrogen demand, renewable production, and electricity spot market prices. 

This study assessed the effectiveness of a number of commonly used planning methods for decision-making in the context of a Hydrogen Purchase Agreement (HPA) offtake contract with a single industrial customer. A market-focused capacity expansion model is used for optimising simultaneously equipment capacities and electricity hedging strategy, and the results were assessed with respect to the levelised cost of hydrogen (LCOH).

Due to the high amount of uncertainty in day-ahead market prices, planning methods incorporating risk aversion were shown to be essential for ensuring a level of stability in LCOH. 

A common practice of limiting the resale of hedging energy on the spot market in modelling was shown to be problematic when choosing between hedging options, as it biases hedging options with a lower resale value.

For a Fixed Demand HPA contract, deterministic planning methods combined with intuitive risk-averse heuristic decisions can be effective. Stochastic methods were found to be necessary to obtain good planning decisions when the HPA contract conditions permit the offtaker to vary their consumption profile, or when access to green subsidies is possible (with the added complexity of their time-matching requirements).

Demand uncertainty in HPA contracts was also seen to necessitate increased equipment sizing (of the electrolyser and hydrogen storage). This increased equipment sizing may be a serious obstacle for some project developers facing difficulties in obtaining financing for the initial investment costs, or in obtaining sufficient hydrogen storage on-site. As such, in the negotiation phase of an HPA, both consumer and producer should consider what forms and levels of demand uncertainty are acceptable to both parties, with the knowledge that less strict contracts may necessitate increased sale prices to fund increased equipment CAPEX costs.

Oversizing electrolyser capacity was not found to be necessary if a Fixed Demand HPA contract is used and the project developer is able to obtain a technologically and geographically diverse set of PPA's with complementary production uncertainties. 

In summary, project developers should be wary when using decisions obtained from deterministic or rule-based modelling to analyse the viability of a particular project. As seen in the case study in this paper, failure to use stochastic modelling can result in decisions that perform over 30\% worse during scenario stress-testing in some contexts. In short, viable projects may be discarded simply because they have been analysed using ineffective planning decisions that perform poorly when tested in adverse conditions.

The validity of the results could be improved by including the presence of short-term production and price forecast uncertainties. Incorporating additional electrolyser dynamics could also have an impact on procurement strategies (potentially favouring less uncertain hedging instruments like futures). Uncertainty in futures price risk premiums could also be considered, as well as the sensitivity of the results to changes in the cost of capital. 

Future work could consider including offtake contract conditions within the optimisation problem, and consider other methods of sharing risk between the consumer and producer, such as price indexing. Use of this model in a multi-stage formulation would also allow decisions to benefit from the evolution of uncertainties over the project life and implement a staged increase of electrolyser capacity.\\

\section{Acknowledgements} \label{sec:acknowledgements}
This study was carried out as part of a collaboration between Mines Paris (PSL University) and Verso Energy. The authors thank the Verso Energy Offers and Partnerships team for providing consultative support, especially Arthur Auxenfants, Jean-Baptiste Martin, and Julien Guiet. The problem is implemented using proprietary modelling software Versys, Verso Energy. Many thanks also to the reviewers of this article who provided many detailed and useful recommendations.

\section{Data Availability} \label{sec:data_availability}
All scenarios used for in-sample and out-of-sample testing are available at \url{https://github.com/owen-verso/Hedging-Hydrogen-Data-Package}.

\appendix

%\appendix

%\include{9D_Data_Availability}


\bibliographystyle{elsarticle-harv.bst} 
\bibliography{references_V5}

\end{document}